\documentclass{article}
\usepackage{amsmath,stmaryrd,amssymb,amsfonts}

\usepackage[hyperindex, colorlinks=true]{hyperref}
\catcode`\@=11
\@addtoreset{equation}{section}

\newtheorem{Theorem}{Theorem}[section]
\newtheorem{Proposition}{Proposition}[section]
\newtheorem{Lemma}{Lemma}[section]
\newtheorem{Corollary}{Corollary}[section]
\newtheorem{Definition}{Definition}[section]
\newtheorem{Remark}{Remark}[section]
\newcommand{\bTheorem}[2]{\bigskip \begin{Theorem}[#2] \label{T#1}}
\newcommand{\eT}{\end{Theorem} \bigskip }
\newcommand{\bProposition}[1]{\bigskip \begin{Proposition} \label{P#1}}
\newcommand{\eP}{\end{Proposition} \bigskip }
\newcommand{\bLemma}[1]{\bigskip \begin{Lemma} \label{L#1}}
\newcommand{\eL}{\end{Lemma} \bigskip }
\newcommand{\bCorollary}[1]{\bigskip \begin{Corollary} \label{C#1}}
\newcommand{\eC}{\end{Corollary} \bigskip }
\newcommand{\bDefinition}[1]{\bigskip \begin{Definition} \label{D#1}}
\newcommand{\eD}{\end{Definition} \bigskip}
\newcommand{\bRemark}[1]{\bigskip \begin{Remark} \label{R#1}}
\newcommand{\eR}{\end{Remark} \bigskip}
\newcommand{\bFormula}[1]{\begin{equation} \label{#1}}
\newcommand{\eF}{\end{equation}}

\newcommand{\bProof}{{\bf Proof: }}
\newcommand{\qed}{\hfill $\boxempty $ \bigskip }

\newcommand{\Ov}[1]{ \overline{ #1 } }

\newcommand{\vc}[1]{ {\bf #1} }

\newcommand{\Div}{{\rm div}_x}
\newcommand{\Curl}{{\bf curl}_x}
\newcommand{\vr}{\varrho}
\newcommand{\vu}{\vc{u}}

\newcommand{\D}{{\cal D}}
\newcommand{\vp}{\varphi}
\newcommand{\dx}{{\rm d}x}
\newcommand{\dt}{{\rm d}t}

\newcommand{\intO}[1]{ \int_{\Omega} #1 \ \dx }

\newcommand{\Grad}{\nabla_x}

\newcommand{\tab}{\,\,\,}

\newcommand{\Om}{\Omega}
\newcommand{\weak}{\rightharpoonup}
\newcommand{\vrho}{\varrho}
\newcommand{\mcw}{\mathcal{W}}

\newcommand{\norm}[1]{\left\Vert#1\right\Vert}

\newcommand{\Set}[1]{\left\{#1\right\}}

\newcommand{\Dom}{\ensuremath{(0,T)\times \Om}}
\newcommand{\Top}[1]{\operatorname{\mathcal{T}}_{h}\left[#1 \right]}

\def\RR{\mathbb{R}}
\def\pat{\partial_t}
\def\nax{\nabla_x}
\def\nav{\nabla_\xi}
\def\eps{\varepsilon}

\title{\bf On the Dynamics of a Fluid-Particle Interaction Model:
The Bubbling Regime}

\author{J. A. Carrillo\thanks{ICREA - Departament de Matem\`atiques, Universitat
Aut\`onoma de Barcelona, 08193 Bellaterra, Spain. E-mail:
carrillo@mat.uab.es}, T. Karper\thanks{Department of Mathematical Sciences,
Norwegian University of Science and Technology, N-7491 Trondheim, Norway.
E-mail: karper@math.ntnu.no }, and K. Trivisa\thanks{Department of Mathematics,
University of Maryland, College Park, MD 20742-4015, U.S.A.
E-mail: trivisa@math.umd.edu }}

\begin{document}

\maketitle


\begin{abstract}
This article deals with the issues of global-in-time existence and
asymptotic analysis of a fluid-particle interaction model in the
so-called bubbling regime. The mixture occupies the physical space
$\Omega \subset \mathbb{R}^3$ which may be unbounded. The system
under investigation describes the evolution of particles dispersed
in a viscous compressible fluid and is expressed by the
conservation of fluid mass, the balance of momentum and the
balance of particle density often referred as the Smoluchowski
equation. The coupling between the dispersed and dense phases
is obtained through the drag forces that the fluid and the
particles exert mutually by action-reaction principle. We show that
solutions exist globally in time under reasonable physical
assumptions on the initial data, the physical domain, and the
external potential. Furthermore, we prove the large-time
stabilization of the system towards a unique stationary state fully
determined by the masses of the initial  density of particles
and fluid and the external potential.
\end{abstract}
\bigskip
\setcounter{tocdepth}{2}
\section{Introduction}
\label{S1}

Fluid-particle interaction arise in many practical applications in
biotechnology, medi\-ci\-ne \cite{BBJM}, reactive gases formed by
fuel droplets in combustion \cite{A,AOr,W,W2}, recycling and
mineral processes \cite{BeBuKaTo03,SG}, and atmospheric pollution
\cite{V}. Aerosols and sprays can be modelled by fluid-particle
type interactions in which bubbles of suspended substances are
seen as solid particles. Two-phase flow hydrodynamic models have
 also been proposed \cite{BDM}.

Here, we focus on a particular system derived by formal
asymptotics from a mesoscopic description. This is based on a
kinetic equation for the particle distribution of Fokker-Planck
type coupled to fluid equations. Different macroscopic equations
can be obtained as scaling limits, see \cite{CG06} for a complete
description of the modelling issues. In these models, the fluid is
either incompressible \cite{GJV1,GJV2} or compressible
\cite{CG06}. The coupling between the kinetic and the fluid
equations is obtained through the friction forces that the fluid
and the particles exert mutually. The friction force is assumed to
follow Stokes law and thus is proportional to the relative
velocity vector, i.e., is proportional to the fluctuations of the
microscopic velocity $\xi\in\RR^3$ around the fluid velocity field
$\vu$. More precisely, the cloud of particles is described by its
distribution function $f_\eps(t, x,\xi)$ on phase space,
which is the solution to the 
dimensionless Vlasov-Fokker-Planck equation
\begin{equation}\label{feps}
\pat f_\eps + \frac{1}{\sqrt\eps}\Big(\xi\cdot \nax f_\eps
-\nabla_x\Phi\cdot\nabla_ \xi f_\eps\Big)= \frac1\eps \ {\mathrm
{div}}_\xi\Big( \big(\xi- \sqrt\eps u_\eps\big)f+ \nav
f_\eps\Big).
\end{equation}
Here, $\epsilon > 0$ is a dimensionless parameter and we have a
 a drag force independent of the fluid density $\vr_\eps,$ but
proportional to the relative velocity of the fluid and the
particles given by $\xi-u_\eps(t,x)$. The right hand-side
of the moment equation in the Navier-Stokes system takes into
account the action of the cloud of particles on the fluid through
the forcing term
$$
F_\eps = \int_{\RR^3}
\left(\frac{\xi}{\sqrt\eps}-u_\eps(t,x)\right) f(t,x,\xi) \,d\xi.
$$
The density of the particles $\eta_\eps(t,x)$  is
related to the probability distribution function $f_\eps(t,
x,\xi)$ through the relation
$$
\eta_\eps(t,x) = \int_{\mathbb{R}^3} f_\eps(t,x,\xi) \, d\xi.
$$

The well-posedness of this kinetic-hyperbolic coupled
fluid-particle system has been addressed in \cite{MV1} in the case
of compressible models for the fluid equations. There are two
different scaling limits for this model, the so-called: bubbling
and flowing regimes. They correspond to the diffusive
approximation of the kinetic equation, bubbling regime, written in
\eqref{feps}, and the strong drag force and strong Brownian motion
for the flowing regime. This last regime has been studied in
\cite{MV2} where it is obtained rigorously as the limit from the
mesoscopic description for local in time solutions and initial
data bounded away from zero for the densities. We also refer to
\cite{CGL} for asymptotic preserving numerical schemes in relation to these scaling limits.

In all the above mentioned studies, the viscosity of the fluid was neglected
although it is the source of the drag forces. The viscosity is
present in the dimensionless systems although negligible as noted
in \cite[Remark 3]{CG06}. In this work, we will deal with the
resulting formal macroscopic fluid-particle system obtained
through the scaling limit in \eqref{feps} as $\eps\to 0$ by the
standard Hilbert-expansion procedure. In this scaling limit,
particles are supposed to have a negligible density with respect
to the fluid, and thus, due to buoyancy effects, they will
typically move upwards in a system under gravity, from where it
bears the name of "bubbling". This situation is typically
complemented with no-flux boundary conditions in a bounded domain.
More generally, we can ask ourselves under what conditions on the
external potential in unbounded domains we can assert convergence
towards stationary integrable states. We refer to \cite{CG06,CGL}
for a detail account of the physical meaning and validity of the
scaling limits.

Summarizing, the state of such flows in this macroscopic
description is, in general, characterized by the variables: the
total mass density $\vr(t,x)$, the velocity field $
\vu(t,x)$, as well as the density of particles in the mixture
$\eta(t,x)$, depending on the time $t \in (0,T)$ and the
Eulerian spatial coordinate $x \in \Omega \subset \RR^3.$

In this section, we present the primitive conservation equations
governing fluid-particle flows in the bubbling regime. These
equations  express the conservation of  mass, the balance of
momentum,  and the balance of particle densities often referred as
the {\em Smoluchowski equation}:
\bFormula{i1}
\partial_t \vr + \Div (\vr \vu) = 0,
\eF
\vspace{0.01in}
\bFormula{i2}
\partial_t (\vr \vu) + \Div (\vr \vu \otimes \vu) + \Grad (p(\vr) + \eta)  -\mu \Delta \vc{u} - \lambda \Grad \Div \vu
= -(\eta + \beta \vr) \Grad \Phi,
 \eF
 \vspace{0.01in}
 \bFormula{i3}
\partial_t \eta + \Div(\eta(\vu - \Grad \Phi)) - \Delta \eta =0.
\eF

\noindent Here, $p$ denotes the pressure $p(\vr) = a
\vr^{\gamma}$, $a > 0$, $\gamma > 1$, $\beta \neq 0$, and $\Phi$
denotes the external potential (typically incorporating gravity
and boyancy). 

In this paper, we require the potential to satisfy
suitable confinement conditions ({\bf HC}) (see Section 2), which does not
limit the physical relevance of our results. The viscosity
parameters $\mu>0$ and $\lambda +\frac{2}{3} \mu\ge 0$ are
nonnegative constants and $\beta > 0$ if $\Om$ is unbounded.

Another macroscopic effect is that the total pressure function in
the momentum equation depends on both the particle and the fluid
densities $p(\vr) + \eta$. We impose the no-slip boundary
condition for the velocity vector leading to no-flux for the fluid
density through the boundaries and the no-flux condition for the
particle density \bFormula{i4} \vu|_{\partial \Omega} = \Grad \eta
\cdot \nu + \eta \Grad \Phi \cdot \nu =0 \tab \mbox{on} \tab
(0,T)\times
\partial\Omega, \eF with $\nu$ denoting the outer normal vector to
the boundary $\partial \Omega.$ Our problem is supplemented with
the initial data $\{\vr_0, \vc{m}_0, \eta_0\}$ such that
 \bFormula{i5}
\begin{split}
     \vr(0,x) & = \vr_0 \in L^{\gamma}(\Omega)\cap L_+^1(\Om), \\
     (\vr \vu)(0,x) &= \vc{m}_0 \in L^{\frac{6}{5}}(\Omega) \cap L^1(\Om),\\
          \eta(0,x) &= \eta_0 \in L^2(\Omega)\cap L_+^1(\Om).
\end{split}
 \eF
Motivated by the stability arguments in \cite{CG06}, the numerical
investigation presented in \cite{CGL}, as well as a number of
studies on numerical experiments and scale analysis on the
proposed model (see \cite{B04}), we investigate the issues of
global existence and asymptotic analysis for fluid-particle
interaction flows providing a rigorous mathematical theory based
on the principle of balance laws. The total energy of the mixture
is given by
\bFormula{i6} E(\eta, \vr, \vu)(t):=
\!\intO{\!\left[\frac{1}{2} \vr(t) |\vu(t)|^2 + \frac{a}{\gamma
-1} \vr^{\gamma}(t) + (\eta \log \eta)(t) + (\beta \vr + \eta)(t)
\Phi\right]\!\!} \eF
At the formal level, the total energy can be
viewed as a Lyapunov function satisfying the {\em energy
inequality} \bFormula{i7} \frac{dE}{dt} + \intO{\Big[\mu |\Grad
\vu|^2 + \lambda |\Div \vu|^2 + |2\Grad  \sqrt{\eta} +
\sqrt{\eta}\Grad \Phi|^2\Big]} \le 0. \eF Therefore, it is
reasonable to anticipate that, at least for some sequences $t_n
\to \infty,$
$$
\eta(t_n) \to \eta_s, \tab \vr(t_n) \to \vr_s, \tab \vr \vu (t_n) \to 0,
$$
where $\eta_s, \vr_s$ satisfy the stationary problem
$$
\Grad(p(\vr_s) + \eta_s) = - (\eta_s + \beta \vr_s) \Grad \Phi
\tab \mbox{on} \tab \Omega.
$$
The energy estimate written in the form
$$
E(t)  + \int_0^T (\| \Grad \vu \|^2_{L^2} + \|\Div \vu\|^2_{L^2}) dt
+ \int_0^T \intO{ |2\Grad \sqrt{\eta} + \sqrt{\eta}\Grad \Phi|^2} dt  \le E(0),
$$
now implies that \bFormula{i9} |2\Grad  \sqrt{\eta_s} +
\sqrt{\eta_s}\Grad \Phi|^2 =0. \eF The aim of this paper is to
show that, in fact, any weak solution converges to a fixed
stationary state as time goes to infinity, more precisely,
$$ \vr(t) \to \vr_s \tab \mbox{strongly in} \tab L^{\gamma}(\Omega),$$
$$ \mbox{ess}\sup_{\tau > t} \intO{\vr(\tau) |\vu(\tau)|^2} \to 0, $$
$$\eta(t) \to \eta_s \tab \mbox{strongly in} \tab  L^p(\Omega), $$
as $t \to \infty$ under the confinement  hypothesis on the domain
$\Omega$ and the external potential $\Phi$ given in {\bf (HC)}
(cf. Definition \ref{Dconfinement}). 

Indeed, it can be shown that the sequences $(\vr_n, \eta_n)$ of the  time shifts defined as
$$\vr_n(t,x) := \vr(t+ \tau_n,x), \,\,\, \tau_n \to \infty,$$
$$\eta_n(t,x) := \eta(t+ \tau_n,x), \,\,\, \tau_n \to \infty,$$
contain subsequences, denoted by the same index w.l.o.g., such
that
$$\vr_n \to \vr_s \tab \mbox{strongly in } \tab L^1_{loc}((0,1) \times \Omega),$$
and
$$\eta_n \to \eta_s \tab \mbox{strongly in } \tab L^{p_1}((0,T) ; L^{p_2}) ( \Omega)  \tab \mbox{for some} \tab p_1, p_2 >1, $$
where $(\vr_s, \eta_s)$ solve the stationary problem
\bFormula{i10} \left\{
\begin{array}{l}
\Grad p(\vr_s) = - \beta \vr_s \Grad \Phi,\\
\tab \tab \,\, \Grad \eta_s = - \eta_s \Grad \Phi.
\end{array}
\right. \eF It is worth noting that the confinement hypothesis is
both necessary and sufficient for the stationary problem
\eqref{i10} to admit a unique solution $(\vr_s, \eta_s)$ (Section \ref{S4}).

The main ingredients of our approach can be formulated  as follows:
\vspace{-0.1in}

\begin{itemize}
\item A suitable {\it weak formulation} of the underlying physical
principles is introduced. Physically grounded hypotheses are
imposed on the system: The mixture occupies the physical space
$\Omega \subset \mathbb{R}^3.$ The boundary conditions are chosen
in such a way so that  the dissipation of energy is guaranteed,
whereas  the pressure of the mixture takes into account both the
density of the fluid and the density of particles.

\vspace{-0.1in} \item {\it A priori} estimates are established
based solely on the boundedness of the initial energy and entropy
of the system.

\vspace{-0.1in} \item A suitable approximating scheme is
introduced for the construction of  the solution based on  a two
level approximating procedure: The first level involves an {\em
artificial pressure approximation}, whereas the second level
approximation employs a {\em time discretization scheme}. The
sequence of approximate solutions is constructed with the aid of a
fixed point argument coupling the time-discretized compressible
isentropic Navier-Stokes equations to a discretization in 
time of the equation for $\eta$.

\vspace{-0.1in} \item Physically grounded hypotheses are imposed on
the domain $\Omega$ and the external potential  $\Phi$
(confinement hypotheses {\bf (HC)}). The analysis in the present
article treats both  the case of a {\em bounded} physical domain
$\Omega$ as well as the case of  an  {\em unbounded}  domain.   We
remark that  the confinement hypothesis {\bf (HC)} on $(\Omega,
\Phi)$ plays a crucial role in providing control of the negative
contribution of the physical entropy $\eta \log \eta$ in the
free-energy bounds for unbounded domains.

\vspace{-0.1in} \item High integrability properties for the
density need to be established for the limit passage in the family
of approximate solutions and in particular in taking the vanishing
artificial pressure limit. We remark  that in the present context,
the potential $ \Phi$ is not integrable on  unbounded domains. To
deal with this new difficulty  we employ the {\em Fourier
multipliers} in the spirit of \cite{EF70,LI4}, while taking into
consideration the new features of our problem.

\vspace{-0.1in} \item We remark that both the total fluid mass and
the total particle mass  are constants of motion. In particular,
we are able to conserve the total masses also in the large-time
limit allowing us to uniquely determine the long time asymptotics (cf~\cite{FP9}).

\end{itemize}

The paper is organized as follows. In Section \ref{S2} we collect
all the necessary hypotheses imposed on the external potential
(confinement hypotheses {\bf (HC)}),  we   present the notion of
{\em free energy solutions} and the main results of this article.
Section \ref{S3} is devoted to the proof of  the global existence
of weak solutions (Theorem \ref{T1}). First,  a suitable
approximation scheme based on  an {\em artificial pressure
approximation} and on  a {\em time discretization scheme}  is
introduced.  The remaining section is
devoted to  the limit passage in the family of approximate
solutions. The most delicate part of the analysis concerned with
the vanishing artificial pressure limit is presented in Subsection
\ref{S3.5}. The large time asymptotic analysis is described in
Section \ref{S4}.

\section{Free-energy solutions and main results}\label{S2}

In this work, we analyse the existence and large-time asymptotics
of certain type of weak solutions to the two-phase flow problem
(\ref{i1})-(\ref{i3}) coupled with no-flux boundary conditions
(\ref{i4}) in two different geometrical constraints of interest in
the applications: for bounded domains and for unbounded domains
under confinement conditions due to the external potential. We
will collect all assumptions concerning the geometry $\Omega$ and
the external potential $\Phi$ under the generic name of {\em
confinement conditions}. Let us remark that the external potential
$\Phi$ is always defined up to a constant. Therefore, for bounded
from below external potentials $\Phi$, we can always assume
without loss of generality, by adding a suitable constant, that
\bFormula{Hinf} \inf_{x \in \Omega} \Phi(x) = 0 . \eF

\bDefinition{confinement} Given a domain $\Omega\in C^{2,\nu}$,
$\nu>0$, $\Omega\subset\RR^3$, and given a bounded from below
external potential $\Phi:\Omega \longrightarrow \RR^+_0$
satisfying \eqref{Hinf}, we will say that $(\Omega,\Phi)$ verifies
the confinement hypotheses {\bf (HC)} for the two-phase flow
system (\ref{i1})-(\ref{i3}) coupled with no-flux boundary
conditions (\ref{i4}) whenever:
\begin{itemize}
\item[] {\bf (HC-Bounded)} If $\Omega$ is bounded, $\Phi$ is
bounded and Lipschitz continuous in $\bar{\Omega}$ and the
sub-level sets $[\Phi < k]$ are connected in $\Omega$ for any
$k>0$.

\item[] {\bf (HC-Unbounded)} If $\Omega$ is unbounded, we assume
that $\Phi \in W^{1,\infty}_{loc}(\Omega)$, $\beta> 0$, the
sub-level sets $[\Phi < k]$ are connected in $\Omega$ for any
$k>0$,
\begin{equation*}
e^{-\Phi/2}\in L^1 (\Omega),
\end{equation*}
and
\begin{equation}\label{potentialreq}
|\Delta \Phi(x)| \leq c_1|\Grad \Phi(x)| \leq c_2\Phi(x), ~  |x|
> R,
\end{equation}
for some large $R>0$.
\end{itemize}
\eD \bRemark{} The confinement assumption {\bf (HC)} has physical
relevance in our setting as it is verified for several domains
$\Om$ with $\Phi$ being the gravitational potential. For instance,
\begin{enumerate}
	\item when
    $\Omega = \{x \in  \RR^3 \,|~(x_1,x_2) \in [a,b]^2,\,  x_3 \in [0,H]\} ~\text{ and } ~\Phi(x) =  g x_3,$
    where $\beta = 1 - \frac{\vr_F}{\vr_P}$.
    
    \item when
    $\Omega = \{x \in  \RR^3 \,|~(x_1,x_2) \in [a,b]^2,\, x_3 > 0\} ~\text{ and } ~\Phi(x) =  g x_3,$
    where $\beta = 1 - \frac{\vr_F}{\vr_P}$ and $\vr_F < \vr_P$.
    \item when $\Om = \mathbb{R}^3 \setminus \overline{B(0,R)}$ and $\Phi(x) = g|x|$, where $B(0,R)$ is the
    ball centered at the origin with radius $R$ and $\beta > 0$.
\end{enumerate}
Here, $\vr_F$ and $\vr_P$ are the typical mass density of fluid and particles, respectively. Remark 
that  1. corresponds to the standard bubbling case (see \cite{CG06}) in
which particles move upwards due to buoyancy.
\eR

Let us now specify the type of weak solutions for the two-phase
flow system (\ref{i1})-(\ref{i3}) we will be dealing with.

\label{S2} \bDefinition{D2a} Let us assume that $(\Omega,\Phi)$
satisfy the confinement hypotheses {\bf (HC)}, we say that $\{
\varrho, {\bf u}, \eta\}$ is a free-energy solution of problem
(\ref{i1})-(\ref{i3}) supplemented with boundary data for which
(\ref{i4}) holds and initial data $\{ \vr_0, \vc{m}_0, \eta_0 \}$
satisfying (\ref{i5}) provided that the following hold:
\begin{itemize}
\item $\vr \ge 0$ represents a renormalized solution of equation
(\ref{i1}) on   $(0,\infty) \times \Omega$:
 For any test function $\vp \in \D([0,T) \times
\Ov{\Omega})$, any $T>0$, and any $b$ such that
\[
b \in L^{\infty} \cap C [0, \infty),\ B(\vr) = B(1) + \int_1^{\vr}
{b(z) \over z^2} \ {\rm d}z,
\]
the following integral identity holds:
\end{itemize}
\vspace{-0.3cm} \bFormula{ii1} \int_0^\infty \intO{ \Big(   B(\vr)
\partial_t \vp +  B(\vr) \vu \cdot \Grad \vp - b(\vr) \Div \vu \vp
\Big)} \ \dt = - \intO{ B(\vr_0) \vp (0, \cdot) } \eF
\begin{itemize}
\item The balance of momentum holds in distributional  sense,
namely \bFormula{ii2} \int_0^\infty \intO{ \Big( \vr \vu \cdot
\partial_t \vc{v} + \vr \vu \otimes \vu : \Grad \vc{v} + (p(\vr) +
\eta) \ \Div \vc{v} \Big) } \ \dt = \eF
\[
 \int_0^\infty \intO{ (\mu
 \Grad \vu \Grad \vc{v} + \lambda \Div \vu \Div  \vc{v}  -
 (\eta + \beta \vr) \Grad \Phi \cdot \vc{v}} \ \dt - \intO{ \vc{m}_0 \cdot \vc{v} (0, \cdot)\!},
\]
for any test function $\vc{v} \in \D([0,T); \D (\Ov{\Omega};
R^3))$ and any $T>0$ satisfying $\vec \vp |_{\partial \Omega} =
0.$

All quantities appearing in (\ref{ii2}) are supposed to be at
least integrable. In particular, the velocity field $\vu$ belongs
to the space $L^2(0,T; W^{1,2}(\Omega; R^3))$, therefore it is
legitimate to require $\vc{u}$ to satisfy the boundary conditions
(\ref{i4}) in the sense of traces.

\item $\eta \geq 0$ is a weak solution of \eqref{i3}. That is, the integral identity \bFormula{ii3} \int_0^\infty\!\!\!
\intO{\eta \partial_t \vp + \eta \vu \cdot \Grad \vp - \eta \Grad
\Phi \cdot \Grad \vp - \Grad \eta \Grad \vp} dt =- \intO{ \eta_0
\vp (0, \cdot)} \eF is satisfied for test functions $\vp \in
{\mathcal D}([0,T)\times \bar{\Omega})$ and any $T>0$.

All  quantities  appear in (\ref{ii3}) must be at least integrable
on $(0,T) \times \Omega$. In particular, $\eta$ belongs
 $L^2([0,T]; L^3(\Omega)) \cap L^1(0,T; W^{1,\frac{3}{2}}(\Omega)).$

\item Given the total free-energy of the system by
$$
E(\vr, \vu, \eta)(t) := \intO{\left( \frac{1}{2} \vr |\vu|^2 +
\frac{a}{\gamma -1} \vr^{\gamma} + \eta \log \eta + (\beta \vr +
\eta) \Phi\right)},
$$
then $E(\vr, \vu, \eta)(t)$ is finite and bounded by the initial
energy of the system, i.e., $E(\vr, \vu, \eta)(t)\leq E(\vr_0,
\vu_0, \eta_0)$ a.e. $t>0$. Moreover, the following free
energy-dissipation inequality holds
\end{itemize}
\eD \vspace{-0.5cm}
\begin{equation}\label{ii4}
\int_0^\infty\!\!\!\intO{\left( \mu |\Grad \vu|^2 + \lambda |\Div
\vu|^2 + |2\Grad \sqrt{\eta} + \sqrt{\eta} \Grad
\Phi|^2\right)}~dt \leq E(\vr_0, \vu_0, \eta_0).
\end{equation}

\vspace{0.5cm}

\noindent Now, we have all ingredients to state the main results
of this work.

\bTheorem{1}{Global Existence} Let us assume that $(\Omega,\Phi)$
satisfy the confinement hypotheses {\bf (HC)}. Then, the problem
(\ref{i1})-(\ref{i3}) supplemented with boundary conditions
(\ref{i4})  and  initial data  satisfying (\ref{i5}) admits a weak
solution $\{ \vr , \vu , \eta\}$ on $(0,\infty   ) \times \Omega$
in the sense of Definition \ref{DD2a}. In addition,
\begin{itemize}
    \item[i)] the total fluid mass and particle mass given by
    $$
    M_\vr(t) = \int_\Om \vr(t, \cdot)~dx \qquad \mbox{and} \qquad
    M_\eta (t) = \int_\Om \eta(t, \cdot)~dx,
    $$
    respectively, are constants of motion.

    \item[ii)] the density satisfies the higher integrability result
    $$
    \vrho \in L^{\gamma + \Theta}((0,T) \times \Om), \text{ for any $T>0$},
    $$
    where $\Theta = \min \{\frac{2}{3}\gamma - 1, \frac{1}{4}\}$.
\end{itemize}
\eT

We can now completely characterize the large time behavior of
free-energy solutions to (\ref{i1})-(\ref{i6}).

\bTheorem{2}{Large-time Asymptotics} Let us assume that
$(\Omega,\Phi)$ satisfy the confinement hypotheses {\bf (HC)}.
Then, for any free-energy solution $(\vr, \vu, \eta)$ of the
problem (\ref{i1})-(\ref{i3}), in the sense of Definition
\ref{DD2a}, there exist universal stationary states
$\vr_s(x)$, $\eta_s(x),$ such that
\[
\left\{
\begin{array}{l}
\vr(t)\rightarrow \vr_s   ~\mbox{strongly  in}~~L^{\gamma}(\Omega),\\ \\
 \displaystyle \mathrm{ess}\sup_{\tau > t} \int_{\Omega} \vr(\tau) |\vu(\tau)|^2 \, dx \to 0, \\ \\
\eta(t)\rightarrow\eta_s ~\mbox{strongly  in}~~L^{p_2}(\Omega) \,\, \mbox{for} \,\, p_2 > 1,\\ \\
\end{array}
\right.
\]
as $t\rightarrow\infty$, where $(\eta_s,\vr_s)$ are characterized
as the unique free-energy solution of the stationary state
pro\-blem: \bFormula{i9} \left\{
\begin{array}{l}
\Grad p(\vr_s) = - \beta \vr_s \Grad \Phi,\\[2mm]
\tab \tab \,\, \Grad \eta_s = - \eta_s \Grad \Phi,
\end{array}
\quad \int_\Om \vr_s~dx = \int_\Om \vr_0~dx, \quad \int_\Om
\eta_s~dx = \int_\Om \eta_0~dx. \right. \eF given by the formulas
\begin{equation*}
\vr_s = \left(\frac{\gamma-1}{a\gamma}\left[-\beta \Phi +
C_\vr)\right]^+\right)^\frac{1}{\gamma -1 }, \qquad \eta_s =
C_\eta\exp(-\Phi),
\end{equation*}
where $C_\eta$ and $C_\vr$ are uniquely given by the initial
masses due to \eqref{i9}. \eT


\section{Global-in-time existence}\label{S3}

This section is devoted to the proof of the   existence result
(Theorem \ref{T1})  which will be achieved by patching
local-in-time solutions with the aid of suitable global estimates.

\begin{Lemma}\label{lemma:existence}
The conclusion of Theorem \ref{T1} holds true on any time-space
cylinder $[0,T)\times \Om$, where $T>0$ is any given finite time.
\end{Lemma}

Taking Lemma \ref{lemma:existence} for granted, a weak solution
verifying Theorem \ref{T1} can be constructed as follows.

\noindent \textbf{Proof of Theorem \ref{T1}:} Fix any $T_0>0$.
From Lemma \ref{lemma:existence}, we have the existence of a weak
solution $(\vr_1, \vu_1, \eta_1)$ on  $[0,2T_0)$. To proceed, we
will need the function
\begin{equation*}
    \zeta^\epsilon(t) =
    \begin{cases}
        1, & t \in [0,T_0-\epsilon], \\
        \operatorname{dist}(t, T_0), & t \in (T_0-\epsilon, T_0), \\
        0, & \text{otherwise}.
    \end{cases}
\end{equation*}
By setting $\vp\zeta^\epsilon$ as a test-function in the
continuity formulation \eqref{ii1}, we obtain
\begin{equation*}
    \begin{split}
    &-\int_{T_0-\epsilon}^{T_0}\int_\Om B(\vr_1)\vp \epsilon^{-1} ~dxdt \\
    &\qquad \qquad + \int_0^{T_0} \intO{ \Big(\zeta^\epsilon B(\vr_1)
    \partial_t \vp  +  \zeta^\epsilon B(\vr_1) \vu_1 \cdot \Grad \vp - \zeta^\epsilon b(\vr_1) \Div \vu_1 \vp
    \Big)} \ \dt  \\
    & = - \intO{ B(\vr_0) \vp (0, \cdot)\zeta^\epsilon(0) }, \quad \forall \vp \in C_c^\infty([0,T]\times \overline{\Om}).
    \end{split}
\end{equation*}
Sending $\epsilon \rightarrow 0$ in this formulation yields
\begin{equation*}
    \begin{split}
        &\int_0^{T_0} \intO{ \Big( B(\vr_1)
        \partial_t \vp  +   B(\vr_1) \vu_1 \cdot \Grad \vp -  b(\vr_1) \Div \vu_1 \vp
        \Big)} \ \dt  \\
        &\qquad \qquad = \int_\Om B(\vr(T_0,\cdot))\vp(T_0, \cdot) ~dxdt - \intO{ B(\vr_0) \vp (0, \cdot)}, \quad \forall \vr \in C_c^\infty([0,T]\times \overline{\Om}).
    \end{split}
\end{equation*}
Hence, $(\vr_1, \vu_1)$ is a renormalized solution of the
continuity equation on the closed interval $[0,T_0]$ with
additional boundary terms at $t=T_0$. By applying similar
arguments to the momentum equation \eqref{ii2} and particle
density equation \eqref{ii3}, we conclude that $(\vr_1, \vu_1,
\eta_1)$ is a weak solution on $[0,T_0]$ (with additional boundary
terms at $t= T_0$).

By the various uniform in time bounds (cf.~\eqref{ii4}),
$(\vr_1,\vu_1, \eta_1)$ at time $T_0$ has sufficient integrability
to serve as initial data for a new solution $(\vr_2, \vu_2,
\eta_2)$ defined on $[T_0, 2T_0] \times \Om$. The energy of the
new solution is less than the initial energy $E(\vr_0, \vu_0,
\eta_0)$ and the triple $(\vr, \vu, \eta)$ given by
$$
(\vr, \vu, \eta)(t) = (\vr_i, \vu_i, \eta_i), \quad t \in
((i-1)T_0, iT_0],
$$
is a weak solution on $[0,2T_0]\times \Om$. A solution for all
times is  readily obtained by iterating this process. \qed

In the remaining parts of this section we prove Lemma
\ref{lemma:existence}. The overall strategy can be outlined as
follows: First, we prove that there exists a sequence of
approximate solutions $\{\vrho_{\delta,h}, \vc{u}_{\delta,h},
\eta_{\delta,h} \}_{h>0}$ (see the ensuing subsection for
details). Then, by first sending $h \rightarrow 0$ and
subsequently $\delta \rightarrow 0$ in this sequence, we prove
that in the limit we obtain a weak solution to the system
\eqref{i1}--\eqref{i3} in the sense of Definition \ref{DD2a}.

\subsection{Approximation scheme}
The approximation scheme is realized in two layers. In the first
layer we add a regularization term to the pressure. 

\begin{Definition}[Artificial pressure approximation]\label{def:level1}
For a given $\delta > 0$, we say that the triple $\{\vr_{\delta},
\vc{u}_{\delta}, \eta_{\delta}\}$ is a weak solution to the
artificial pressure approximation scheme provided that
$\{\vr_{\delta}, \vc{u}_{\delta}, \eta_{\delta}\}$ is a weak
solution in the sense of Definition \ref{DD2a} but with the
pressure $p_{\delta}(\vrho_{\delta})$ given by
$$
    p_{\delta}(\vrho_{\delta}) = p(\vrho_{\delta}) + \delta \vrho_{\delta}^6,
$$
and initial data $(\vr_{\delta,0}, m_{0,\delta}, \eta_{\delta, 0})
= \left(\vr_0 , m_0, \eta_0\right)\star \kappa_\delta,$
 with $\kappa_\delta$ denoting the standard smoothing kernel.
\end{Definition}

The second layer of approximation is a discretization of the
equations in time. To define this layer of approximation, we shall
make use of the space
$$
\mathcal{W}(\Om) = L^1(\Om)\cap L^6(\Om)\times
W^{1,2}_{0}(\Om)\times W^{1,\frac{3}{2}}(\Om)\cap L^3(\Om)\cap
L^1(\Om).
$$

\begin{Definition}[Time discretization scheme]\label{def:level2}
Let $\delta > 0$ be fixed. Given a time step $h> 0$, we discretize
the time interval $[0,T]$ in terms of the points $t^k = kh$, $k=0,
\ldots, M$, where we assume that $Mh = T$. Now, we sequentially
determine functions
$$
\{\vrho^k_{\delta,h},\vc{u}_{\delta,h}^k, \eta_{\delta,h}^k\} \in
\mathcal{W}(\Om),\quad k=1, \ldots, M,
$$
such that:
\begin{itemize}
\item{} The time discretized continuity equation,
\begin{equation}\label{eq:timedisc-vrho}
d_{t}^h[\vrho^k_{\delta,h}] + \Div (\vrho_{\delta,h}^k
\vc{u}_{\delta,h}^k) = 0,
\end{equation}
holds in the sense of distributions on $\overline{\Om}$. \item{}
The time discretized momentum equation with artificial pressure,
\begin{equation}\label{eq:timedisc-momentum}
    \begin{split}
        &d_{t}^h[\vrho^k_{\delta,h} \vc{u}^k_{\delta,h}] + \Div (\vrho_{\delta,h}^k\vc{u}^k_{\delta,h}\otimes \vc{u}^k_{\delta,h}) - \mu \Delta \vc{u}^k_{\delta,h} - \lambda \Grad \Div \vc{u}^k_{\delta,h} \\
        &\qquad + \Grad \left(p_{\delta}(\vrho_{\delta,h}^k) + \eta^k_{\delta,h}\right) = - (\beta \vrho_{\delta,h}^k + \eta_{\delta,h}^k)\Grad \Phi,
    \end{split}
\end{equation}
holds in the sense of distributions on $\Om$. \item{} The time
discretized particle density equation,
\begin{equation}\label{eq:timedisc-eta}
d_{t}^h[\eta_{\delta,h}^k] + \Div \left(\eta_{\delta,h}^k
(\vc{u}_{\delta,h}^k - \Grad \Phi) \right) - \Delta
\eta_{\delta,h}^k = 0,
\end{equation}
holds in the sense of distributions on $\overline{\Om}$.
\end{itemize}
In the above equations, $d_{t}^h[\phi^k] = \frac{\phi^{k} -
\phi^{k-1}}{h}$ denotes implicit time stepping.
\end{Definition}

For each fixed $h>0$, the time discretized solution
$\{\vrho_{\delta,h}^k, \vc{u}_{\delta,h}^k,
\eta_{\delta,h}^k\}_{k=1}^M$ is extended to the whole of
$(0,T)\times \Om$ by setting
\begin{equation}\label{eq:level2-ext}
(\vrho_{\delta, h}, \vc{u}_{\delta, h}, \eta_{\delta, h})(t) =
(\vrho_{\delta,h}^k, \vc{u}^k_{\delta,h}, \eta_{\delta,h}^k),
\qquad t \in (t^{k-1},t^{k}],\ k=1, \ldots, M.
\end{equation}
In addition, we set the initial data
$$
\left(\vrho_{\delta,h}(0), (\vrho_{\delta,0}\vc{u}_{\delta,h})(0),
\eta_{\delta, h}(0)\right) = (\vrho_{\delta,0}, m_{\delta,0},
\eta_{\delta,0}).
$$

Next two subsections treat the existence of a well-defined
approximation scheme in both relevant cases $\Omega$ bounded or
not under the confinement conditions {\bf (HC)}.


\subsection{Well-defined Approximations: $\Om$ bounded}

To prove existence of a solution to the time discretized
approximate equations
\eqref{eq:timedisc-vrho}--\eqref{eq:timedisc-eta} we will utilize
a fixed point argument. For this purpose we define  the operator
$$
    \Top{\cdot}: \mathcal{W}(\Om) \rightarrow \mathcal{W}(\Om),
$$
as the solution $( \vrho, \vc{u}, \eta) = \Top{\psi, \vc{z},
\zeta}$ to the system of equations:
\begin{equation}\label{eq:fix-timedisc-vrho}
\begin{split}
\frac{\vrho- \vrho^{k-1}}{h} + \Div (\vrho\vc{u}) & = 0,
\end{split}
\end{equation}
\begin{equation}\label{eq:fix-timedisc-momentum}
\begin{split}
&\frac{\varrho\vc{u} - \vr^{k-1} \vu^{k-1} }{h}+ \Div (\varrho\vc{u}\otimes \vc{u}) -\mu \Delta \vc{u} - \lambda \Grad \Div \vc{u} \\
&\qquad = - \Grad(p_{\delta}(\varrho)  + \zeta) - (\beta \varrho+
\zeta)\Grad \Phi,
\end{split}
\end{equation}
and
\begin{equation}\label{eq:fix-timedisc-eta}
\begin{split}
 \frac{\eta- \eta^{k-1}}{h} + \Div\left(\eta (\vc{z} - \Grad \Phi)\right) - \Delta \eta &= 0,
 \end{split}
\end{equation}
in the sense of distributions on $\Om$, where $(\vr^{k-1},
\vu^{k-1}, \eta^{k-1}) \in \mcw(\Om)$ is given.

Observe that for a fixed $h> 0$, $\delta > 0$, and $k=1,\ldots,
M$, any fixed point of the operator $\Top{\cdot}$ will be a
distributional solution to the time discretized approximate
equations \eqref{eq:timedisc-vrho}--\eqref{eq:timedisc-eta}. i.e
the fixed point
$$
(\vrho_{\delta, h}^k, \vc{u}_{\delta, h}^k, \eta_{\delta, h}^k) =
\Top{\vrho_{\delta, h}^k, \vc{u}_{\delta, h}^k, \eta_{\delta,
h}^k}
$$
is precisely the desired solution at each time step.

\subsubsection{$\mathcal{T}_h$ is well-defined}
The existence of functions $(\vrho, \vc{u},\eta) \in \mcw(\Om)$
satisfying $( \vrho, \vc{u},\eta) = \Top{ \psi, \vc{z},\zeta}$
follows from Lemmas \ref{lemma:stationary} and \ref{lem:etaexistence} below.
\begin{Lemma}[\cite{LI4}, Theorem 6.1]\label{lemma:stationary}
Let $h>0$ be fixed arbitrary and assume that
$\frac{1}{h}\vrho^{k-1} \in L^\gamma(\Om)\cap L^6(\Om)$,
 $\frac{1}{h}m^{k-1} \in L^\frac{6}{5}(\Om)\cap L^p(\Om)$ for some $p> 1$, $ \zeta \in W^{1, \frac{3}{2}}(\Om)$,
and $\Grad \Phi \in L^p(\Om)$, for all $p$. There exists a weak
solution $(\vrho, \vc{u}) \in L^\gamma(\Om)\cap L^6(\Om)\times
W^{1,2}_{0}(\Om)$ of
\eqref{eq:fix-timedisc-vrho}--\eqref{eq:fix-timedisc-momentum}
satisfying
\begin{equation}\label{eq:pressL2}
p_{\delta}(\vrho) \in L^{2}(\Om),
\end{equation}
and
\begin{equation}\label{eq:effviscstrongconv}
\Grad \Curl \vc{u},\  \Grad \left(\Div \vc{u}- \frac{1}{\lambda +
\mu}p_{\delta}(\vrho) \right), \ \in L^q(K),
\end{equation}
where $q = \frac{15}{11} > \frac{6}{5}$ and $K \subset \Omega$
is any compact subset.
\end{Lemma}

{\it Notation.-} Throughout the paper we use overbars to denote
weak limits, just to illustrate this notational use, we will often
use convexity arguments based on the following classical lemma
\cite{EF70}:

\begin{Lemma}\label{lem:prelim}
Let $O$ be a bounded open subset of $\RR^M$ with $M\ge 1$. Suppose
$g\colon \RR\rightarrow (-\infty,\infty]$ is a lower semicontinuous
convex function and $\Set{v_n}_{n\ge 1}$ is a sequence of
functions on $O$ for which $v_n\weak v$ in $L^1(O)$, $g(v_n)\in
L^1(O)$ for each $n$, $g(v_n)\weak \overline{g(v)}$ in $L^1(O)$.
Then $g(v)\le \overline{g(v)}$ a.e.~on $O$, $g(v)\in L^1(O)$, and
$\int_O g(v)\ dy \le \liminf_{n\to\infty} \int_O g(v_n) \ dy$. If,
in addition, $g$ is strictly convex on an open interval
$(a,b)\subset \RR$ and $g(v)=\overline{g(v)}$ a.e.~on $O$, then,
passing to a subsequence if necessary, $v_n(y)\to v(y)$ for
a.e.~$y\in \Set{y\in O\mid v(y)\in (a,b)}$.
\end{Lemma}
\begin{Lemma}\label{lem:etaexistence}
Assume that $\Om$ is bounded.
 Let $\eta^{k-1} \in L^1(\Om)\cap W^{1, \frac{3}{2}}(\Om)\cap\{\eta^{k-1}\geq 0\}$, and $\vc{z}Ê\in \vc{W}^{1,2}_{0}(\Om)$ be given functions.
Then, for each fixed $h > 0$, there exists a  non-negative
function $\eta \in W^{1,\frac{3}{2}}(\Om)\cap L^1(\Om)$ satisfying
\eqref{eq:fix-timedisc-eta} in the sense of distributions on
$\Omega$. Moreover,
\begin{equation}\label{eq:lem42}
    \int_\Om h^{-1}[\eta \log \eta + \eta\Phi] + |2\Grad \sqrt{\eta} + \sqrt{\eta}\Grad \Phi|^2 ~dx \leq C,
\end{equation}
where the constant $C>0$ depends on $\|\vc{z}\|_{\vc{W}^{1,
2}(\Om)}$, $\|\eta^{k-1}\|_{L^3(\Om)}$, and
$\|\Phi\|_{W^{1,\infty}(\Om)}$.
\end{Lemma}
\bProof For each  $\epsilon > 0$, we let $\vc{z}^\epsilon =
\vc{z}\star \kappa_\epsilon$, where $\star$ is the convolution
product and $\kappa_\epsilon$ is the standard smoothing kernel.
From \cite[Proposition {\bf 4.29}]{NS1}, we can assert the
existence of a unique weak solution $\eta^\epsilon \in
W^{2,2}(\Om)$, $\eta^\epsilon \geq 0$, to the linear elliptic
equation
\begin{equation}\label{eq:ep1}
\frac{\eta^\epsilon- \eta^{k-1}}{h} + \Div\left(\eta^\epsilon
(\vc{z}^\epsilon - \Grad \Phi)\right) - \Delta \eta^\epsilon = 0
\text{ in }\Om,
\end{equation}
satisfying the boundary condition
\begin{equation*}
    \eta^\epsilon \Grad \Phi\cdot \nu + \Grad \eta^\epsilon\cdot \nu = 0.
\end{equation*}
By integrating over $\Om$ (using the boundary condition), we
observe that
\begin{equation*}
    \int_\Om \eta^\epsilon~dx = \int_\Om \eta^{k-1}~dx \leq C,
\end{equation*}
and hence $\eta^\epsilon \in L^1(\Om)$ independently of
$\epsilon$.

Let the sequence $\{B_l\}_{l=1}^\infty$ be given by
\begin{equation*}
    B_l(y) =
    \begin{cases}
        \log y, & y > l^{-1}, \\
        \log l^{-1}, & y\leq l^{-1}.
    \end{cases}
\end{equation*}
Moreover, let $\Om_l = \{x \in \Om: \eta^\epsilon (x) > l^{-1}\}$
and $\Om^c_l = \Om \setminus \Om_l$. Multiply \eqref{eq:ep1} with
$ B_l(\eta^\epsilon)$ and integrate by parts to obtain
\begin{equation}\label{eq:seherja}
    \begin{split}
        &\int_{\Om_l} h^{-1}\left[\eta^\epsilon \log \eta^\epsilon - \eta^{k-1}\log \eta^\epsilon\right]~dx
        - \int_{\Om^c_l}h^{-1}\left[\eta^\epsilon \log l  - \eta^{k-1}\log l\right]~dx\\
        &\qquad
        + \int_{\Om_l} \Grad \Phi \cdot \Grad \eta^\epsilon + \frac{1}{\eta^\epsilon}|\Grad \eta^\epsilon|^2
        - \vc{z}^\epsilon \Grad\eta^\epsilon ~dx
         = 0.
    \end{split}
\end{equation}
Applying the identity $\tfrac{1}{\eta^\epsilon}|\Grad
\eta^\epsilon|^2 = 4|\Grad \sqrt{\eta^\epsilon}|^2$ and reordering
terms in \eqref{eq:seherja}
\begin{equation*}
    \begin{split}
        &\int_{\Om_l} h^{-1}\eta^\epsilon \log \eta^\epsilon
         + 4|\Grad \sqrt{\eta^\epsilon}|^2~dx \\
        &\qquad =
        \int_{\Om_l^c}h^{-1}\left[\eta^\epsilon \log l  - \eta^{k-1}\log l\right]~dx
        + \int_{\Om_l}h^{-1}\eta^{k-1}\log \eta^\epsilon~dx\\
        &\qquad \qquad - \int_{\Om_l} 2\sqrt{\eta^\epsilon}\Grad \Phi \cdot \Grad \sqrt{\eta^\epsilon}
        - 2\sqrt{\eta^\epsilon}\vc{z}^\epsilon \cdot \Grad\sqrt{\eta^\epsilon} ~dx \\
        &\qquad \leq \int_{\Om_l^c}h^{-1} \log l^\frac{1}{l}~dx + \int_{\Om_l}h^{-1}\eta^{k-1}\eta^\epsilon~dx \\
        &\qquad \qquad + 2\|\Grad \sqrt{\eta^\epsilon}\|_{L^2(\Om_l)}^2\left(\|\sqrt{\eta^\epsilon}\Grad \Phi\|_{L^2(\Om_l)}
         + \|\sqrt{\eta^\epsilon}\vc{z}^\epsilon\|_{L^2(\Om_l)}\right) \\
        &\qquad \leq h^{-1}\log l^\frac{1}{l} |\Om_l^c| + \int_{\Om_l}h^{-1}\eta^{k-1}\eta^\epsilon~dx
                + 2\|\sqrt{\eta^\epsilon}\|_{L^2(\Om_l)}^2 \\
        &\qquad \qquad
                    + \|\sqrt{\eta^\epsilon}\vc{z}^\epsilon\|_{L^2(\Om_l)}^2 + \|\sqrt{\eta^\epsilon}\Grad \Phi\|_{L^2(\Om_l)},
    \end{split}
\end{equation*}
where the last inequality is the Cauchy inequality. The last term
is bounded since $\eta^\epsilon, \eta^{k-1} \in L^1(\Om)$ and
$\Phi \in W^{1,\infty}(\Om)$, independently of $\epsilon$ and $l$.
Using the Sobolev embedding $\vc{W}^{1,2}(\Om) \subset
\vc{L}^6(\Om)$ together with $\eta^{k-1} \in L^3(\Om)$, we achieve
\begin{equation}\label{eq:gohome}
    \begin{split}
        &\int_{\Om_l} h^{-1}\eta^\epsilon \log \eta^\epsilon
        + 2|\Grad \sqrt{\eta^\epsilon}|^2~dx \\
        &\qquad  \leq C + \int_{\Om_l} \eta^\epsilon\left(h^{-1}\eta^{k-1} + |\vc{z}^\epsilon|^2\right)~dx
        \leq  C(1 + \|\eta^\epsilon\|_{L^\frac{3}{2}(\Om_l)}).
    \end{split}
\end{equation}
Applying the H\"older inequality, Sobolev embedding, and Young's
inequality (with epsilon), respectively,  we
 bound the last term as follows
\begin{equation*}
    \begin{split}
        \|\eta^\epsilon\|_{L^\frac{3}{2}(\Om_l)}
        \leq \|\eta^\epsilon\|_{L^1(\Om)}^\frac{1}{3}\|\sqrt{\eta^\epsilon}\|_{L^4(\Om_l)}^\frac{4}{3}
        \leq C\|\Grad \sqrt{\eta^\epsilon}\|_{L^2(\Om_l)}^\frac{4}{3}
        \leq C + \beta \|\Grad \sqrt{\eta^\epsilon}\|_{L^2(\Om_l)}^2.
    \end{split}
\end{equation*}
Inserting this expression in \eqref{eq:gohome} and fixing $\beta$
small, we gather
\begin{equation}\label{eq:particle1}
    \int_{\Om_l} h^{-1}\eta^\epsilon \log \eta^\epsilon  + |\Grad \sqrt{\eta^\epsilon}|^2~dx
    \leq C,
\end{equation}
where the constant $C$ is independent of both $\epsilon$ and $l$.
Since $\Omega$ is bounded, it follows that
\begin{equation*}
    \eta^\epsilon\log \eta^\epsilon \in L^1(\Om), \quad \eta^\epsilon \in W^{1, \frac{3}{2}}(\Om),
\end{equation*}
independently of $\epsilon$.

Since $W^{1,\frac{3}{2}}(\Om)$ is compactly embedded in $L^p(\Om)$
for all $p<3$, we can conclude that, passing to a subsequence if
necessary,
\begin{equation*}
    \eta^\epsilon \weak \eta ~\text{ in }~W^{1,\frac{3}{2}}(\Om), \qquad \eta^\epsilon \rightarrow \eta ~\text{ in }~L^p(\Om),~p<3,
\end{equation*}
as $\epsilon \rightarrow 0$. Consequently, we can take the limit
$\epsilon \rightarrow 0$ in \eqref{eq:ep1} to conclude that
$\eta$ satisfies \eqref{eq:fix-timedisc-eta} in the sense of
distributions on $\Omega$.

Next, since $\eta^\epsilon \rightarrow \eta$ in $L^p(\Om)$, $p<3$,
we can also conclude that
$$
\sqrt{\eta^\epsilon} \rightarrow \sqrt{\eta} ~\text{ in
}~L^q(\Om),~ q < 6,
 \quad \Grad \sqrt{\eta^\epsilon} \weak \Grad \sqrt{\eta}~\text{ in }~L^2(\Om),
$$
as $\epsilon \rightarrow 0$. Using this, we can send $\epsilon
\rightarrow 0$ in \eqref{eq:particle1} to obtain
\begin{equation*}\label{eq:particle2}
    \int_\Om h^{-1}\eta \log \eta   + \overline{|\Grad \sqrt{\eta}|^2}~dx
    \leq C.
\end{equation*}
By weak lower semi-continuity $\overline{\|\Grad
\sqrt{\eta}\|^2_{L^2(\Om)}} \geq \|\Grad
\sqrt{\eta}\|^2_{L^2(\Om)}$. Thus,
\begin{equation*}\label{eq:particle2}
    \int_\Om h^{-1}\eta \log \eta
     + |\Grad \sqrt{\eta}|^2
    \leq C.
\end{equation*}
Then, using the identity
\begin{equation*}
    2\Grad \Phi \cdot \Grad \eta^\epsilon + \eta^\epsilon |\Grad \Phi|^2 + 4|\Grad \sqrt{\eta^\epsilon}|^2
    = |2\Grad \sqrt{\eta^\epsilon} + \sqrt{\eta^\epsilon}\Grad \Phi|^2,
\end{equation*}
we are led to the conclusion
\begin{equation*}
    \begin{split}
        &\int_\Om h^{-1}\left[\eta \log \eta + \eta\Phi \right]
        + |2\Grad \sqrt{\eta^\epsilon} + \sqrt{\eta^\epsilon}\Grad \Phi|^2~dx \\
        & \qquad \leq C + \int_\Om h^{-1}\eta\Phi+  4 \Grad \sqrt{\eta} \cdot \left(\sqrt{\eta}\Grad \Phi \right)
        + \eta|\Grad \Phi|^2~dx \\
        &\qquad \leq C\left(1 + \|\eta\|_{L^1(\Om)}\left(\|\Phi\|_{L^\infty(\Om)}
                    + \|\Grad \Phi\|_{L^\infty(\Om)}^2 + \|\Grad \Phi\|_{L^\infty(\Om)}\|\Grad \sqrt{\eta}\|_{L^2(\Om)}\right)\right) \\
        &\qquad \leq \tilde{C},
    \end{split}
\end{equation*}
which is \eqref{eq:lem42} and the proof is complete.\qed

\vspace{-0.3cm}

\begin{Remark}\label{rem:rem}
    Since $\Grad \sqrt{\eta} \in L^2(\Om)$, the set of particle vacuum regions
    have measure zero; $\left|\{x \in \Om; \eta = 0\} \right| = 0$.
\end{Remark}


\subsubsection{$\mathcal{T}_h$ admits a fixed point}
We now prove that $\Top{\cdot}$ admits a fixed point and
consequently that the time discretization scheme in Definition
\ref{def:level2} is well defined. The key observation made is that
the $L^2$ bound on the pressure enables us to obtain an energy
equality. This equality in turn yields compactness of the operator
$\Top{\cdot}$.

\begin{Lemma}\label{lem:fixed-point}
Assume the case of bounded domain $\Om$. Let $\{\vrho^{k-1},
m^{k-1}, \eta^ {k-1}\} \in L^2(\Om)\times L^\frac{6}{5 }(\Om)\cap
L^p(\Om)\times L^\gamma(\Om), p>1,$ be given functions. Then, for
each fixed $h >0$ there exists a fixed point $\{\vrho, \vc{u},
\eta\} \in \mcw(\Om)$ for the operator $\Top{\cdot}$; i.e
$$
    (\vrho, \vc{u},\eta) = \Top{\vrho, \vc{u},\eta}.
$$
As a consequence, the time discretization scheme given by
Definition \ref{def:level2} is well-defined for bounded domains
$\Om$.
\end{Lemma}

\bProof We will prove the existence of a fixed point by verifying
the postulates of the Schauder corollary to the Schaefer fixed
point theorem \cite{Deimling}; \emph{If the operator $\Top{\cdot}$
is continuous, compact, and the set $\Set{s \in [0,1], x \in
\mcw(\Om): x = \Top{sx}}$ is uniformly bounded, then the operator
$\Top{\cdot}$ admits a fixed point.}

First, we observe that the operator $\Top{\cdot}$ is clearly
bounded and continuous.

Next, we prove that the operator $\Top{\cdot}$ is compact. For
this purpose, let $\Set{0, \vc{z}_{n}, \eta_{n}}_{n=1}^\infty$ be
a sequence such that $(\vc{z}_{n},\eta_{n}) \in_{b}
W^{1,2}_{0}(\Om)\times W^{1,\frac{3}{2}}(\Om)$ for all $n=1,
\ldots, \infty$, and construct a sequence $\Set{\vrho_{n},
\vc{u}_{n}, \eta_{n}}_{n=1}^\infty$ by setting
$$
(\vrho_{n}, \vc{u}_{n}, \eta_{n}) = \Top{0, \vc{z}_{n},
\zeta_{n}}, \quad n  = 1, \ldots, \infty.
$$
Then, from the previous lemmas we know that $(\vr_{n},\vc{u}_{n},
\eta_{n}) \in \mcw(\Om)$ independently of $n$. Hence, we have the
existence of functions $(\vr, \vc{u}, \eta) \in \mcw(\Om)$ such
that, by passing along a subsequence if necessary,
\begin{equation*}
(\vrho_{n}, \vc{u}_{n}, \eta_n) \weak (\vrho, \vc{u}, \eta), \quad
\textrm{in }\mcw(\Om).
\end{equation*}
Moreover, by compact Sobolev embedding, we clearly have the
existence of $(0,\vc{z}, \eta) \in \mcw(\Om)$ such that $\zeta_{n}
\rightarrow \zeta$ a.e in $\Om$ and $\vc{z}_{n} \rightarrow
\vc{z}$ a.e in $\Omega$.

Now, we claim that in fact $(\vrho_{n}, \vc{u}_{n}, \eta_n)
\rightarrow (\vrho, \vc{u}, \eta)$ in $\mathcal{W}(\Om)$, where
$(\vrho, \vc{u}, \eta) = \Top{0, \vc{z}, \zeta}$, and consequently
that the operator $\Top{\cdot}$ is compact.

To prove this claim, we first note that compactness of $\eta_{n}$
in $W^{1,\frac{3}{2}}(\Om)$ is immediate from the linearity of
\eqref{eq:timedisc-eta}. That is, since $\eta_{n} \in W^{1,
\frac{3}{2}}(\Om)$ we have by Sobolev embedding that $\eta_{n}
\rightarrow \eta$ a.e in $\Omega$ and thus by setting $\log \eta -
\log \eta_{n}$ (cf.~Remark \ref{rem:rem}) as test--function in
\eqref{eq:fix-timedisc-eta}  one discovers
\begin{equation*}
\begin{split}
\lim_{n \rightarrow \infty}\int_{\Omega}|\Grad \sqrt{\eta_{n}}|^2
-|\Grad \sqrt{\eta}|^2\ dx = 0,
\end{split}
\end{equation*}
which immediately implies compactness in $W^{1,\frac{3}{2}}(\Om)$.
Moreover, in the limit we have that $\eta$ satisfies
\begin{equation}\label{eq:fix-etadone}
 \frac{\eta- \eta^{k-1}}{h} + \Div\left(\eta (\vc{z} - \Grad \Phi)\right) - \Delta \eta = 0,
\end{equation}
in the sense of distributions on $\overline{\Om}$.

We continue proving compactness of the operator $\Top{\cdot}$ by
first proving strong convergence of the density $\vrho_{n}
\rightarrow \vrho$ a.e in $\Om$. Compactness of the velocity
$\vc{u}_{n}$ in $W^{1,2}_{0}(\Om)$
 will then follow from an energy equality.
In order to prove strong convergence of the density we will need
weak sequential continuity of the effective viscous flux. That is,
\begin{equation}\label{eq:discviscconv}
\begin{split}
&\lim_{n \rightarrow \infty}\int_{\Omega}\left[(\lambda + \mu)\Div \vc{u}_{n} - p_{\delta}(\vrho_{n})\right] \psi \vrho_{n} \ dx \\
&= \int_{\Omega}\left[(\lambda + \mu)\Div \vc{u} -
\overline{p_{\delta}(\vrho)} \right] \psi \vrho \ dx, \quad
\forall \psi \in C_{c}^\infty(\Om).
\end{split}
\end{equation}
Here, \eqref{eq:discviscconv} is an immediate consequence of
\eqref{eq:effviscstrongconv} and
 compact  embedding of Sobolev spaces.

Before we proceed to prove strong convergence of $\vrho_{n}$, we
first note that there is no problem with taking the limit in
\eqref{eq:timedisc-vrho} to obtain in the limit
\begin{equation}\label{eq:fix-vrhodone}
\frac{1}{h}\vrho + \Div (\vc{u} \vrho) = \frac{1}{h}\vrho_{k-1},
\end{equation}
in the sense of distributions on $\overline{\Omega}$. Hence, since
in particular $\vrho \in L^6(\Om)$ and $\vc{u} \in
W^{1,2}_{0}(\Om)$ we can conclude that, for any $B \in
C[0,\infty)\cap C^1(0,\infty)$,
\begin{equation}\label{eq:timedisc-cont-renorm}
\frac{1}{h}\vrho B'(\vrho) + \Div (B(\vrho)\vc{u}) + ((\vrho
B'(\vrho) - B(\vrho)) \Div \vc{u}) =
\frac{1}{h}\vrho_{k-1}B'(\vrho),
\end{equation}
 in the sense of distributions on $\overline{\Om}$.
Then, by setting $B(z) = z \log z$ we obtain the identity
\begin{equation*}
\begin{split}
&\int_{\Omega}\vrho_{n}\log \vrho_{n}- \vrho\log \vrho +  \vrho_{n} - \vrho\ dx \\
&= \int_{\Omega}\vrho \Div \vc{u} - \vrho_{n}\Div \vc{u}_{n} +
\vrho_{k-1}(\log \vrho_{n} - \log \vrho)\ dx.
\end{split}
\end{equation*}
Thus, by taking the limit $n \rightarrow \infty$, we see that
\begin{equation}\label{eq:timedisc-closetocompact}
\begin{split}
&\int_{\Omega}\overline{\vrho \log \vrho} - \vrho \log \vrho \, dx \\
&= \lim_{l \rightarrow \infty}\lim_{n \rightarrow
\infty}\int_{\Omega}\psi_{l}\vrho \,\Div \vc{u} -
\psi_{l}\vrho_{n}\Div \vc{u}_{n}Ê\ dx +
\int_{\Om}\vrho_{k-1}(\overline{\log \vrho} - \log \vrho)\ dx,
\end{split}
\end{equation}
where $\psi_{l} \in C_{c}^\infty(\Om)\cap \Set{\psi_{l}\geq 0}$ is
such that $\psi_{l} = 1$ on the set $ \Set{x \in \Om;
\operatorname{dist}(x, \partial \Om) > \frac{1}{l}}$. By an
application of \eqref{eq:discviscconv} we see that for each
$l=1\ldots \infty$,
\begin{equation*}
\begin{split}
&\lim_{n \rightarrow \infty}\int_{\Omega}\psi_{l}\vrho\, \Div
\vc{u} - \psi_{l}\vrho_{n}\Div \vc{u}_{n}Ê\ dx = \lim_{n
\rightarrow \infty}\int_{\Omega}
\psi_{l}\left(p_{\delta}(\vrho_{n})\vrho -
p_{\delta}(\vrho_{n})\vrho_{n} \right) \ dx \leq 0,
\end{split}
\end{equation*}
where the last inequality follows from the convexity of
$p_{\delta}(\vrho_{n})$. Similarly,  the concavity of $\log
\vrho_{n}$ yields
\begin{equation*}
\int_{\Om}\vrho_{k-1}(\overline{\log \vrho} - \log \vrho)\ dx \leq
0.
\end{equation*}
However, then \eqref{eq:timedisc-closetocompact} together with the
convexity of $z\log z$ gives
\begin{equation*}
\begin{split}
0 \leq \int_{\Omega}\overline{\vrho \log \vrho} - \vrho \log \vrho
dx \leq 0.
\end{split}
\end{equation*}
which immediately yields $\vrho_{n} \rightarrow \vrho$ a.e in
$\Om$ due to Lemma \ref{lem:prelim}.

We are now in a position to prove that $\vc{u}_{n} \rightarrow
\vc{u}$ in $W^{1,2}_{0}(\Om)$. For this purpose, we first set
$\vc{u}_{n}$ as test--function in \eqref{eq:timedisc-momentum} to
obtain the identity
\begin{equation}\label{eq:gammlefar}
\begin{split}
&\lim_{n \rightarrow \infty}\int_{\Omega}\frac{\varrho_{n} |\vc{u}_{n}|^2 + \vrho_{k-1}|\vc{u}_{n}|^2}{2h} - \frac{\vc{u}_{n}m_{k-1}}{h} + \mu |\Grad \vc{u}_{n}|^2 + \lambda |\Div \vc{u}_{n}|^2 \ dx\\
&\qquad =\lim_{n \rightarrow \infty} \int_{\Omega}\left(p_{\delta}(\vrho_{n}) +  \zeta_{n}\right) \Div \vc{u}_{n} - (\zeta_{n}  + \beta \vrho_{n})\Grad \Phi\cdot \vc{u}_n \ dx \\
&\qquad = \int_{\Omega}\left(p_{\delta}(\vrho) +  \zeta\right) \Div
\vc{u} - (\zeta  + \beta \vrho)\Grad \Phi\cdot \vc{u} \ dx.
\end{split}
\end{equation}
Note that in addition to the strong convergence of $\eta_{n}$ and
$\vrho_{n}$ we really need \eqref{eq:pressL2} to conclude this.

Now, since $\vrho_{n} \rightarrow \vrho$, $\zeta_{n} \rightarrow
\zeta$, and $\vc{u}_{n} \rightarrow \vc{u}$ a.e in $\Om$ there is
no problem with passing to the limit in
\eqref{eq:timedisc-momentum} to obtain
\begin{equation*}
\begin{split}
\frac{\varrho\vc{u} - m_{k-1} }{h}&+ \Div (\varrho\vc{u}\otimes \vc{u}) -\mu \Delta \vc{u} - \lambda \Grad \Div \vc{u} \\
& = - \Grad (p_{\delta}(\varrho)  + \zeta) - (\beta \varrho+
\zeta)\Grad\Phi,
\end{split}
\end{equation*}
in the sense of distributions on $\Om$. Hence, in view of
\eqref{eq:fix-etadone} and \eqref{eq:fix-vrhodone}, we can
conclude that $(\vrho, \vc{u}, \eta) = \Top{ 0, \vc{z}, \zeta}$.
Moreover, by using $\vc{u}$ as a test--function for this equation
we obtain the identity
\begin{equation*}
\begin{split}
\int_{\Omega} &\left[\frac{\vrho \vc{u}^2+ \vrho_{k-1}\vc{u}^2}{2h} -\frac{\vc{u}m_{k-1}}{h}+ \mu |\Grad \vc{u}|^2 + \lambda |\Div \vc{u}|^2 \right] dx   \\
&= \int_{\Omega}(p_{\delta}(\vrho) + \zeta)\Div \vc{u}- (\zeta + \beta \vrho) \Grad \Phi\cdot \vc{u} \ dx \\
&= \int_{\Omega}\frac{\vrho \overline{\vc{u}^2}+
\vrho_{k-1}\overline{\vc{u}^2}}{2h} -\frac{\vc{u}m_{k-1}}{h}+ \mu
\overline{|\Grad \vc{u}|^2} + \lambda \overline{|\Div \vc{u}|^2}
dx,
\end{split}
\end{equation*}
where the last equality is \eqref{eq:gammlefar}. This can only
happen if $\Grad \vc{u}_{n} \rightarrow \Grad \vc{u}$. Thus, we
can conclude that $\Top{\cdot}$ is a compact operator.

Finally, let $s \in [0,1]$ be arbitrary and assume that there
exists a triple $(\vrho, \vc{u}, \eta) \in \mcw(\Om)$ such that
$(\vrho, \vc{u}, \eta)  = \Top{s \vr,  s \vc{u}, s \eta}$. Then,
by setting $\vc{u}$ as test-function in
\eqref{eq:timedisc-momentum} we get the identity
\begin{equation}\label{eq:fixedpoint-uenergy}
\begin{split}
&\int_{\Omega}\frac{\vrho \vc{u}^2+ \vrho_{k-1}\vc{u}^2}{2h} -\frac{\vc{u}m_{k-1}}{h}+ \mu |\Grad \vc{u}|^2 + \lambda |\Div \vc{u}|^2 dx   \\
&\qquad = \int_{\Omega} (p_{\delta}(\vrho) + s\eta)\Div \vc{u}  -
(s\eta + \beta \vrho) \Grad\Phi\cdot \vc{u} \ dx.
\end{split}
\end{equation}
Using $B(z) = \Pi_{\delta}(z) := \frac{1}{\gamma - 1}p(\vrho) +
\frac{\delta}{5}\vrho^6$ as the renormalization function in
\eqref{eq:timedisc-cont-renorm} we also have the identity
\begin{equation}\label{eq:fixedpoint-vrhoenergy}
    \begin{split}
        -\int_{\Omega}p_{\delta}(\vrho)\Div \vc{u}\ dx
        &= \frac{1}{h}\int_{\Omega}\Pi_{\delta}'(\vrho)(\vrho - \vrho_{k-1})\ dx.
    \end{split}
\end{equation}
Similarly, using $\beta \Phi$ as test-function in
\eqref{eq:fix-timedisc-vrho} gives
\begin{equation*}
\int_{\Omega}\beta \frac{\vrho - \vrho_{k-1}}{h}\Phi - \beta
\vrho\, \vc{u} \cdot \Grad \Phi dx = 0.
\end{equation*}
Using $\Phi + \log \eta$ as a test-function in
\eqref{eq:fix-timedisc-eta} (cf.~Lemma \ref{lem:etaexistence})
gives
\begin{equation}\label{eq:fixedpoint-etaen}
\begin{split}
0 &=\int_{\Omega}\frac{\eta - \eta_{k-1}}{h}(\Phi + \log \eta) -
s\eta \vc{u}\cdot \Grad\Phi
+ s\eta \Div \vc{u} \\
&\qquad + |2\Grad \sqrt{\eta} + \sqrt{\eta}\Grad \Phi|^2\ dx.
\end{split}
\end{equation}
By combining
\eqref{eq:fixedpoint-vrhoenergy}--\eqref{eq:fixedpoint-etaen}, we
obtain the identity
\begin{equation*}
\begin{split}
&\int_{\Omega}p_{\delta}(\vrho) \Div \vc{u} + s\eta \Div \vc{u} - (s\eta + \beta \vrho) \Grad \Phi\cdot \vc{u} \ dx \\
&= \int_{\Omega}-\Pi_{\delta}'(\vrho)\frac{(\vrho -
\vrho_{k-1})}{h}
-\frac{\eta - \eta_{k-1}}{h}(\Phi + \log \eta) \\
&\qquad -\beta \frac{\vrho - \vrho_{k-1}}{h}\Phi -  |2\Grad
\sqrt{\eta} + \sqrt{\eta}\Grad \Phi|^2 dx.
\end{split}
\end{equation*}
Then, inserting this into \eqref{eq:fixedpoint-uenergy} and
reordering terms gives
\begin{equation}\label{eq:timedisc-energyequality}
\begin{split}
&\int_{\Omega}\frac{\vrho \vc{u}^2+ \vrho_{k-1}\vc{u}^2}{2h}
-\frac{\vc{u}m_{k-1}}{h} + \frac{\vrho -
\vrho_{k-1}}{h}\left(\Pi_{\delta}'(\vrho)
+ \beta\Phi\right)  + \frac{\eta - \eta_{k-1}}{h}(\Phi + \log \eta) dx   \\
&\qquad + \int_{\Omega}\mu |\Grad \vc{u}|^2 + \lambda |\Div \vc{u}|^2 \ dx \\
&= -\int_{\Omega}|2\Grad \sqrt{\eta} + \sqrt{\eta}\Grad \Phi|^2 dx
\leq 0.
\end{split}
\end{equation}
Consequently, we can conclude that
$$
\|\eta\|_{W^{1,\frac{3}{2}}(\Om)} + \|\vrho\|_{L^\gamma(\Om)} +
\|\vc{u}\|_{W^{1,2}(\Om)} \leq C,
$$
where the constant $C$ depends on the data $\eta_{k-1},
\vrho_{k-1}$, $m_{k-1}$, and $\Phi$, together with $h$. However,
$C$ does not depend on the parameter $s$.

We can now conclude the proof since we have proved that the
operator $\Top{\cdot}$ is bounded, continuous, compact, and that
the set $\Set{s \in [0,1], x \in \mcw(\Om): x = \Top{sx}}$ is
uniformly bounded. \qed


\subsection{Well-defined Approximations: $\Om$ unbounded}

At this point we have proved that the approximation scheme is
well-defined on bounded domains. Given an unbounded domain
$\Omega$ and an external potential $\Phi$ satisfying the
assumptions {\bf (HC)}, we can always find an increasing sequence
of domains $\Omega_r$, with $r> 0$ such that $\Omega_r$ are
bounded and $(\Omega_r,\Phi)$ satisfies {\bf (HC)} approximating
$\Omega$ in the sense $\cup_{r>0}\Omega_r =\Omega$. Using the
previous subsection, for any $r>0$, there is a solution on
$\Omega_r$. In this subsection, we prove that we can send $r
\rightarrow \infty$ to obtain a solution in $\Omega$. The
following lemmas will be of use in the sequel.

\subsubsection{Consequences of confinement: $\Omega$ unbounded}

We show in this part how to control the negative contribution of
the physical entropy $\eta\log\eta$ in the free-energy bounds for
unbounded domains $\Omega$. Here, the confinement conditions {\bf
(HC)} on $(\Omega,\Phi)$ are crucial. Most of these lemmas can be
seen in \cite{Do99} but we include them here for the sake of
completeness.
 We first start by a classical lemma in kinetic theory.

\begin{Lemma}\label{logminus}
Assume that $(\Omega,\Phi)$ satisfy the hypotheses {\bf (HC)}. For
any density $\eta\in L^1_+(\Omega)$,
\[\label{Ineq:EntrUnif}
\int_{\Omega}\eta(x)\,\log_- \eta(x)\,\, dx \leq \frac12
\int_{\Omega}\Phi(x)\eta(x)\,\ dx + \frac 1e \int_{\Omega}
e^{-\Phi(x)/2} \, dx\;.
\]
\end{Lemma}

\bProof Let $\bar \eta:=\eta\,\chi_{\{\eta\le 1\}}$ and $\bar
M=\int_{\Omega}\bar \eta(x)\, dx\le\int_{\Omega} \eta(x)\, dx=M$.
Then
\[
\int_{\Omega}\bar \eta(x)\,\left(\log\bar \eta(x)+\frac12
\Phi(x)\right)\, dx=\int_{\Omega} [U(x)\log U(x)]\mu\, dx-\bar M
\log Z
\]
where $U:=\bar \eta/\mu$, $\mu(x)= e^{-\Phi(x)/2}/Z$ with
$Z=\int_{\Omega} e^{-\Phi(x)/2} \, d x$. The Jensen inequality
yields
\[
\int_{\Omega} [U(x)\,\log U(x)]\mu \, dx \ge \left(\int_{\Omega}
U(x)\mu \, dx\right)\;\log\left(\int_{\Omega} U(x)\mu \, dx
\right)=\bar M\,\log \bar M\;
\]
and
\begin{align*}
-\int_{\Omega}\eta(x)\,\log_- \eta(x)\,\, dx = \int_{\Omega}\bar
\eta(x)\,\log\bar \eta(x)\, dx & \ge \bar M\log
\bar M - \bar M \log Z -\frac12 \int_{\Omega}\Phi(x)\,\bar \eta(x)\, dx\\
& \ge -\frac Ze-\frac12 \int_{\Omega}\Phi(x)\, \eta(x)\, dx\;,
\end{align*}
from which the desired claim follows. \qed

We can immediately use this previous lemma to conclude the
following consequence.

\begin{Corollary}\label{corentropybound}
Assume that $(\Omega,\Phi)$ satisfy the hypotheses {\bf (HC)}. For
any density $\eta\in L^1_+(\Omega)$, if
$$
\int_{\Omega}\eta(x)\,\log \eta(x)\,\, dx +
\int_{\Omega}\Phi(x)\eta(x)\,\ dx \leq C\, ,
$$
then $\eta \log \eta\in L^1(\Omega)$ and there exists $D>0$
depending on $C$ and $\Phi$ such that
$$
\int_{\Omega}\eta(x)\,\log_+ \eta(x)\,\, dx \leq D \qquad
\mbox{and} \qquad \int_{\Omega}\Phi(x)\eta(x)\,\ dx \leq D\;.
$$
\end{Corollary}

Finally, the above estimates can be used to control the mass of
the densities $\eta$ outside a large ball to avoid loss of mass at
infinity.

\begin{Lemma}\label{lemcontmass}
Given any domain $\Omega$ such that $e^{-\Phi}\in L^1_+(\Omega)$
and any density $\eta\in L^1_+(\Omega)$, then
$$
\int_{\Omega}\eta(x)\,\log \eta(x)\,\, dx +
\int_{\Omega}\Phi(x)\eta(x)\, dx \geq \int_{\Omega}\eta(x)\, dx
\,\log\left(\frac{\int_{\Omega}\eta(x)\,
dx}{\int_{\Omega}e^{-\Phi(x)}\, dx}\right)\, .
$$
As a consequence, if $e^{-\Phi}\in L^1_+(\Omega)$ and
$$
\int_{\Omega}\eta(x)\,\log \eta(x)\,\, dx +
\int_{\Omega}\Phi(x)\eta(x)\,\ dx \leq C\, ,
$$
then, for any $\epsilon>0$ there exists $R>0$ depending on $C$ and
$\Phi$ only such that
$$
\int_{\Omega\cap (\RR^3/B(0,R))}\eta(x)\, dx <\epsilon \,.
$$
\end{Lemma}

\bProof A direct use of Jensen's inequality shows the first
inequality by using the convexity of $x \mapsto x\log x$. To show
the second claim we start by applying the first inequality to the
domain $\Omega_R^c:=\Omega\cap (\RR^3/B(0,R))$ from which we
obtain
\begin{equation}\label{new}
\int_{\Omega_R^c}\eta(x)\, dx
\,\log\left(\frac{\int_{\Omega_R^c}\eta(x)\,
dx}{\int_{\Omega_R^c}e^{-\Phi(x)}\, dx}\right) \leq D
\end{equation}
for some $D>0$, where Lemma \ref{logminus} and Corollary
\ref{corentropybound} were used. Now, we argue by contradiction,
if the second claim were not true, we would have
$$
\exists \epsilon>0 \,\forall R_0>0 \,\exists R>R_0 \mbox{ such
that } \int_{\Omega_R^c}\eta(x)\, dx \geq \epsilon_0.
$$
Since $e^{-\Phi}\in L^1_+(\Omega)$, we can always assume that
$R_0$ is large such that
$$
\int_{\Omega_R^c}e^{-\Phi(x)}\, dx \leq
\int_{\Omega_{R_0}^c}e^{-\Phi(x)}\, dx < \epsilon_0 \leq
\int_{\Omega_R^c}\eta(x)\, dx
$$
and thus due to \eqref{new},
$$
\int_{\Omega_R^c}\eta(x)\, dx \leq \int_{\Omega_R^c}e^{-\Phi(x)}\,
dx \, e^{D/\epsilon_0} \leq \int_{\Omega_{R_0}^c}e^{-\Phi(x)}\, dx
\, e^{D/\epsilon_0} \, .
$$
This leads to a contradiction since the right-hand side can be
made arbitrarily small by taking $R_0$ large enough.\qed

\subsubsection{The approximation scheme is well-defined}
{\it Notation.} - In what follows, we will often obtain a priori
estimates for a sequence $\Set{v_n}_{n\ge 1}$ that we write as
``$v_n\in_{b} X$'' for some functional space $X$. What this really
means is that we have a bound on $\norm{v_n}_X$ that is
independent of $n$.

\begin{Lemma}
Set $(\Om,\Phi)$ satisfying {\bf (HC)} with $\Omega$ unbounded and
let $(\vr_h^{k-1}, \vu_h^{k-1}, \eta^{k-1}_h) \in
\mathcal{W}(\Om)$ be given data. The time discretization scheme
\eqref{eq:timedisc-vrho}-\eqref{eq:timedisc-eta} admits a
distributional solution in the sense of Definition
\ref{def:level2}.
\end{Lemma}
\bProof Let $\{\Omega_r\}_{r>0}$ be an increasing sequence of
domains such that $\cup_{r>0}\Omega_r =\Omega$ and such that, for
each fixed $r$, $\Omega_r$ is bounded and $(\Omega_r,\Phi)$
satisfies {\bf (HC)}. From the results of the previous subsection,
we have the existence of a triplet $(\vr_r, \vu_r, \eta_r)$
satisfying the time discretized equations
\eqref{eq:timedisc-vrho}-\eqref{eq:timedisc-eta} in the sense of
distributions on $\Om_r$. Consequently, we can define a family of
such solutions:
$$
\{\vr_r, \vu_r, \eta_r\}_{R<r<\infty}, \qquad (\vr_r, \vu_r,
\eta_r) \in \mathcal{W}(\Om_r), \text{ for each fixed } r \in
(R,\infty),
$$
where $R$ is fixed according to the requirements on the potential
(see \eqref{potentialreq}). For this construction,
\eqref{eq:timedisc-energyequality} yields
\begin{equation}\label{eq:bs1}
\begin{split}
&\int_{\Om_r}\frac{\vrho_r \vc{u}_r^2+ \vrho^{k-1}\vc{u}_r^2}{2h}
 + \frac{\vrho_r}{h}\left(\Pi_{\delta}'(\vrho_r)
+ \beta\Phi\right)  + \frac{\eta_r}{h}(\Phi + \log \eta_r) dx   \\
&\qquad + \int_{\Om_r}\mu |\Grad \vc{u}_r|^2 + \lambda |\Div \vc{u}_r|^2 + |2\Grad \sqrt{\eta_r} + \sqrt{\eta_r}\Grad \Phi|^2dx \\
&\leq \int_{\Om_r}\frac{\vc{u}_rm_{k-1}}{h} +  h^{-1}\eta_{k-1}\log \eta_r~dx \\
&\leq \epsilon h^{-1} [\|\eta_r\|_{L^2(\Om_r)}^2 +
\|\vu\|_{L^6(\Om_r)}^2] +
      \frac{1}{ 4h\epsilon}\left[\|\eta_{k-1}\|_{L^2(\Om)}^2  + \|m_{k-1}\|_{L^\frac{6}{5}(\Om)}\right].
\end{split}
\end{equation}
In addition, integrating \eqref{eq:timedisc-eta} over $\Om_r$,
\begin{equation*}
    \|\eta_r\|_{L^1(\Om_r)} = \|\eta^{k-1}\|_{L^1(\Om_r)} \leq C,
\end{equation*}
with constant $C$ independent of $r$. An interpolation inequality
and the Sobolev-Poincar\'e inequality allow us to conclude
\begin{equation}\label{eq:bs2}
    \begin{split}
        \|\eta_r\|_{L^2(\Om_r)} &\leq \|\eta_r\|_{L^1(\Om_r)}^\frac{1}{4}\|\eta_r\|_{L^3(\Om_r)}^\frac{3}{4}
        \leq C\|\Grad \sqrt{\eta_r}\|_{L^2(\Om_r)}^\frac{3}{4},
    \end{split}
\end{equation}
where the constant $C$ is independent of $r$. By applying
\eqref{eq:bs2} and the Sobolev-Poincar\'e inequality to
\eqref{eq:bs1}, we conclude
\begin{align}\label{eq:bs3}
        &\int_{\Om_r}\frac{\vrho_r \vc{u}_r^2+ \vrho^{k-1}\vc{u}_r^2}{2h}
         + \frac{\vrho_r}{h}\left(\Pi_{\delta}'(\vrho_r)
        + \beta\Phi\right)  + \frac{\eta_r}{h}(\Phi + \log \eta_r) dx   \\
        &\qquad \!\! + \int_{\Om_r}(\mu - \epsilon) |\Grad \vc{u}_r|^2 + \lambda |\Div \vc{u}_r|^2 + |2\Grad \sqrt{\eta_r} + \sqrt{\eta_r}\Grad \Phi|^2
                -\epsilon|\Grad \sqrt{\eta_r}|^2  dx
        \leq C, \nonumber
\end{align}
where $C$ is independent of $r$.

Next, we observe that \eqref{eq:bs3} yields
\begin{equation*}
    \begin{split}
        &\int_{\Om_r} 4|\Grad \sqrt{\eta_{r}}|^2 + 2\Grad \eta_{r} \Grad \Phi + \eta_{r}|\Grad \Phi|^2~dx
            = \int_{\Om_r} |2\Grad \sqrt{\eta_{r}} + \sqrt{\eta_{r}}\Grad \Phi|^2~dx \leq C.
    \end{split}
\end{equation*}
Reordering terms and integrating by parts
\begin{equation*}
    \begin{split}
         \int_{\Om_r} 4|\Grad \sqrt{\eta_{r}}|^2 + \eta_{r}|\Grad \Phi|^2~dxdt
         \leq C + \int_{\Om_r} \eta_{r} |\Delta \Phi|~dx.
    \end{split}
\end{equation*}
Then, \eqref{potentialreq} and \eqref{eq:bs3} gives
\begin{equation}\label{eq:thecalculation}
    \begin{split}
        &\int_{\Om_R} \eta_{r} |\Delta \Phi|~dx + \int_{\Om_r \setminus \Om_R} \eta_{r} |\Delta \Phi|~dx \\
        &\qquad \leq \|\Delta \Phi\|_{L^\infty(\Om_R)}\int_{\Om_R}\eta_r~Êdx +
        C\int_{\Om_r \setminus \Om_R} \eta_{r} \Phi~dx \leq C,
    \end{split}
\end{equation}
with $C$ independent of $r$.

Setting \eqref{eq:thecalculation} into \eqref{eq:bs3}, fixing
$\epsilon$ small, and applying Corollary \ref{corentropybound}
gives
\begin{equation*}
    \begin{split}
        &\int_{\Om_r}\frac{\vrho_r \vc{u}_r^2+ \vrho^{k-1}\vc{u}_r^2}{2h}
         + \frac{\vrho_r}{h}\left(\Pi_{\delta}'(\vrho_r)
        + \beta\Phi\right)  + \frac{\eta_r}{h}(\Phi + \log_+ \eta_r) dx   \\
        &\qquad + \int_{\Om_r}\mu |\Grad \vc{u}_r|^2 + \lambda |\Div \vc{u}_r|^2 + |\Grad \sqrt{\eta_{r}}|^2 + \eta_{r}|\Grad \Phi|^2 dx
        \leq C.
    \end{split}
\end{equation*}
Since $\|\sqrt{\eta_r}\|_{L^2(\Om_r)}^2= \int_{\Om_r} \eta_{r}~dx
\leq C$,
 the previous estimate allow us to conclude
\begin{equation}\label{eq:etafirstofall}
\eta_{r}  \in_{b} L^3(\Om_r)\cap W^{1,\frac{3}{2}}(\Om_r),
\end{equation}
\begin{align*}
    \vu_r &\in_b L^6(\Om_r)\cap W^{1,2}_0(\Om_r), \qquad
    \vr_r \in_b L^6(\Om_r)\cap L^\gamma(\Om_r).
\end{align*}

Next, let $(\phi, \vc{v}, \psi) \in [C_c^\infty(\Om)]^3$ be
arbitrary and fix a number $\overline{r} \in (R, \infty)$ large
 such that
\begin{equation*}
    \operatorname{supp}(\phi, \vc{v}, \psi) \subset \Om_r.
\end{equation*}
Then, from the previous bounds, we have the existence of functions
$(\vr, \vu, \eta) \in \mathcal{W}(\Om)$ such that, along a
subsequence,
\begin{equation}\label{eq:bs5}
    \begin{split}
        \vr_r \weak \vr, \quad &\text{in }L^6(\Om_{\bar r}), \quad
        \vu_r \weak \vu, \quad \text{in }W^{1,2}(\Om_{\bar r}),  \\
        &\eta_r \weak \eta, \quad \text{in }W^{1,\frac{3}{2}}(\Om_{\bar r})).
    \end{split}
\end{equation}
Moreover, by compact Sobolev  embedding
\begin{equation}\label{eq:bs6}
    \begin{split}
        \eta_r &\rightarrow \eta, \quad \text{in }L^3_{\text{loc}}(\Om_{\bar r}), \qquad
        \vu_r \rightarrow \vu, \quad \text{in } L^p_{\text{loc}}(\Om_{\bar r}), ~p<6.
    \end{split}
\end{equation}
Using \eqref{eq:bs5}, \eqref{eq:bs6}, there is no problem with
sending $r \rightarrow 0$ in
\begin{equation*}
\int_{\Om_{\bar r}}d_t^h[\vrho_{r}] \phi - \vr_r\vc{u}_{r}\Grad
\phi\ dx = 0,
\end{equation*}
\begin{equation*}
\int_{\Om_{\bar r}}d_t^h[\eta_{r}]\psi + (\vc{u}_{r} -
\Grad\Phi)\Grad\psi - \Grad \eta_{r} \Grad \psi\ dxdt = 0,
\end{equation*}
to conclude that $(\vr, \vu, \eta)$ solves
\eqref{eq:timedisc-vrho}-\eqref{eq:timedisc-eta} in the sense of
distributions on $\Om$ (Recall that $(\phi, \psi)$ was chosen
arbitrary). Similarly, we can pass to the limit in the momentum
equation \eqref{eq:timedisc-momentum} to obtain
\begin{equation*}
    \begin{split}
        &\int_{\Om_{\bar r}} d_t^h[\vr\vu]\vc{v} - \vr\vu \otimes \vu :\Grad \vc{v} + \mu \Grad \vu \Grad \vc{v}
        + \lambda \Div \vu  \Div \vc{v}~dx \\
        &\qquad + \int_{\Om_{\bar r}} (\beta \vr + \eta)\Grad \Phi \vc{v}~Êdx
        = \lim_{r \rightarrow \infty} \int_{\Om_{\bar r}} p_\delta(\vr_r)\Div \vc{v}~dx.
    \end{split}
\end{equation*}
Hence, it  remains to prove that $p_\delta(\vr_r) \weak
p_\delta(\vr)$ as $r \rightarrow \infty$. This problem is also the
core problem encountered in the existence analysis of Lions
\cite{LI4}. Due to the regularity properties of $\eta_r$, the
presence of the $\eta$ variable does not impose any additional
difficulties as compared with \cite{LI4}. However, some minor
modifications are needed to incorporate the unbounded potential
$\Phi$. The needed modifications for the present stationary case
are almost identical to those of the non-stationary case. Thus, we
do not give the arguments here but refer the reader to the more
general arguments given in Section 4.5.2. \qed


\subsection{There exists an artificial pressure solution ($h \rightarrow 0$)}
In the previous subsection we proved that, for each $h> 0$ and
$\delta > 0$, we can construct functions $(\vrho_{\delta,h},
\vc{u}_{\delta,h}, \eta_{\delta,h})$  according to  Definition
\ref{def:level2} and \eqref{eq:level2-ext}.  In this section, we
prove that the corresponding sequence $\{(\vrho_{\delta,h},
\vc{u}_{\delta,h}, \eta_{\delta,h})\}_{h>0}$, with $\delta > 0$
fixed, converges as the time step $h \rightarrow 0$  to  an
artificial pressure solution in the sense of Definition
\ref{def:level1}.


\begin{Lemma}
    \label{lemma:simon}
\cite[Corollary 4, p. 85]{Simon}. Let $X \subset B \subset Y$, be
Banach spaces with $X \subset B$ compactly. Then, for $1 \leq p <
\infty$, $\{v: v \in L^p(0,T;X), v_t \in L^1(0,T;Y)\}$ is
compactly embedded in $L^p(0,T;B)$.
\end{Lemma}

The following lemma is a variation of a result due to P.~L.~Lions
(see \cite{Karlsen}).
\begin{Lemma}\label{lemma:dLions}
For a given $T>0$, divide the time interval $(0,T)$ into $M$
points such that $(0,T) = \cup_{k=1}^M(t_{k-1}, t_{k}]$, where
$t_{k} = hk$ and we assume that $Mh = T$. Let
$\{f_{h}\}_{h>0}^\infty$, $\{g_{h}\}_{h>0}^\infty $ be two
sequences such that:
\begin{itemize}
\item{}$\{f_{h}\}_{h>0},\{g_{h}\}_{h>0}$ converge weakly to $f,g$ respectively in $L^{p_{1}}(0,T;L^{q_{1}}(\Om)),$ \\
$L^{p_{2}}(0,T;L^{q_{2}}(\Om))$ where $1 \leq p_{1},q_{1}\leq
\infty$,
$$
\frac{1}{p_{1}} + \frac{1}{p_{2}} = \frac{1}{q_{1}} +
\frac{1}{q_{2}} = 1.
$$
\item{} the mapping $t \rightarrow g_{h}(t,x)$ is constant on each
interval $(t_{k-1}, t_{k}],\ k=1, \ldots, M$. \item{} the discrete
time derivative satisfies,
$$
\frac{g_{h}(t,x) - g_{h}(t-h,x)}{h} \in_{b}
L^1(0,T;W^{-1,1}(\Om)).
$$
\item{} $\{f_{h}\}_{h>0}$ satisfies:
$
\|f_{n} - f_{n}(\cdot + \xi,t)\|_{L^{p_{2}}(0,T;L^{q_{2}}(\Om))}
\rightarrow 0, \textrm{ as }|\xi|\rightarrow 0,
$
uniformly in $n$.
\end{itemize}
Then, $g_{h}f_{h} \weak gf$ in the sense of distributions on
$\Dom$.
\end{Lemma}

\subsubsection{Energy estimates}
Let $\delta > 0$ be fixed and let $\{\vrho_{\delta, h},
\vc{u}_{\delta, h}, \eta_{\delta, h}\}_{h>0}$ be a sequence of
time discretized solutions constructed according to Definition
\ref{def:level2} and \eqref{eq:level2-ext}. Since
\eqref{eq:timedisc-energyequality} holds for every $k$, we can sum
this equality over all $k=1, \ldots, m$, for any $m\in [1,M]$,  to
obtain the energy equality
\begin{equation}\label{eq:disc-energybound}
\begin{split}
&E(\vrho_{\delta, h},\vc{u}_{\delta, h},\eta_{\delta, h})(t^m) \\
&\qquad + \int_{0}^{t^m}Ê\int_{\Omega}\mu |\Grad \vc{u}_{\delta, h}|^2 + \lambda |\Div \vc{u}_{\delta,h}|^2 \ dxdt \\
&\qquad + \int_{0}^{t^m}\int_{\Omega}|2\Grad \sqrt{\eta_{\delta,h}} + \sqrt{\eta_{\delta,h}}\Grad \Phi|^2\ dxdt \\
&\qquad + \frac{1}{\gamma - 1}\sum_{k=1}^m \int_{\Omega}(\gamma - 1)(\vrho_{\delta,h}^k)^\gamma + (\vrho^{k-1}_{\delta,h})^\gamma - \gamma (\vrho_{\delta,h}^{k-1})^{\gamma-1}\vrho_{\delta,h}^k\ dx\\
&\qquad + \frac{\delta}{5}\sum_{k=1}^m \int_{\Omega}5(\vrho_{\delta,h}^k)^6 + (\vrho^{k-1}_{\delta,h})^6 - 6 (\vrho_{\delta,h}^{k-1})^{5}\vrho_{\delta,h}^k\ dx \\
&\qquad + \sum_{k=1}^m \int_{\Om}\eta_{k-1}\log \left(\frac{\eta_{k-1}}{\eta_{k}}\right)\ dx \\
&= E(\vrho_{\delta,0}, \vc{u}_{\delta, 0}, \eta_{\delta,0}),
\end{split}
\end{equation}
where $t^m= mh  \in (0,T)$ and the energy $E(\cdot, \cdot, \cdot)$
is given by \eqref{ii4}. By convexity $(\gamma- 1) a^\gamma +
b^\gamma - \gamma a^{\gamma-1}b \geq 0$ for all $a,b \geq 0,
\gamma > 1$. By concavity of $z \mapsto \log z$,
$$
\int_\Om \eta_{k-1}\log \left(\frac{\eta_{k-1}}{\eta_{k}}\right)
\geq \int_\Om \eta^{k-1} - \eta^{k}~dx = 0.
$$
Due to the confinement conditions, this, Corollary \ref{corentropybound}, and
\eqref{eq:disc-energybound}, allow us to conclude the bounds:
\begin{equation}\label{eq:spaces1}
\begin{split}
\vrho_{\delta,h} |\vc{u}_{\delta,h}|^2 &\in_{b} L^\infty(0,T;L^1(\Om)), \\
\vc{u}_{\delta,h} &\in_{b} L^2(0,T;W^{1,2}_{0}(\Om)), \\
\vrho_{\delta,h} &\in_{b} L^\infty(0,T;L^\gamma(\Om) \cap L^6(\Om)), \\
\eta_{\delta,h} \log \eta_{\delta,h} &\in_{b}
L^\infty(0,T;L^1(\Om)).
\end{split}
\end{equation}

%
By similar arguments as those leading to \eqref{eq:etafirstofall},
we conclude
\begin{equation}\label{eq:etaonceandforall}
\eta_{\delta,h}  \in_{b} L^2(0,T;L^3(\Om))\cap
L^1(0,T;W^{1,\frac{3}{2}}(\Om)),
\end{equation}
independent of both $h$ and $\delta$.

Utilizing the above $h$ independent bounds, it is a simple
exercise to obtain from the system
\eqref{eq:timedisc-vrho}--\eqref{eq:timedisc-eta} the weak time
difference bounds
\begin{equation}\label{eq:timedisc-weaktime}
\begin{split}
d_{t}^h[\vrho_{\delta,h}] &\in_{b} L^\infty(0,T;W^{-1,\frac{3}{2}}(\Om)),\\
d_{t}^h[\eta_{\delta,h}] &\in_{b} L^1(0,T;W^{-1,\frac{3}{2}}(\Om)), \\
d_{t}^h[\vrho_{\delta,h}\vc{u}_{\delta, h}] &\in_{b}
L^1(0,T;W^{-1,1}(\Om)).
\end{split}
\end{equation}

\subsubsection{Convergence}
In view of \eqref{eq:spaces1} and \eqref{eq:etaonceandforall}, we
have the existence of functions $(\vr_\delta, \vu_\delta,
\eta_\delta)$ such that, along a subsequence as $h \rightarrow 0$,
\begin{align*}
    \vr_{\delta, h} &\overset{*}{\weak}\vr_\delta \quad \text{in } L^\infty(0,T;L^\gamma(\Om)), \\
    \vu_{\delta,h} &\weak \vu_\delta \quad \text{in } L^2(0,T;\vc{W}^{1,2}_0(\Om)), \\
    \eta_{\delta,h} &\weak \eta_\delta \quad \text{in } L^2(0,T;L^3(\Om)).
\end{align*}
By virtue of \eqref{eq:timedisc-weaktime}, Lemma
\ref{lemma:dLions} can be applied to yield
\begin{equation}\label{eq:discconv}
\begin{split}
\vrho_{\delta, h}\vc{u}_{\delta,h} \weak
\vrho_{\delta}\vc{u}_{\delta},
 \quad \vrho_{\delta,h}\vc{u}_{\delta,h}\otimes \vc{u}_{\delta,h}\weak \vrho_{\delta}\vc{u}_{\delta}\otimes \vc{u}_{\delta},
\end{split}
\end{equation}
in the sense of distributions on $\Dom$ as $h \rightarrow 0$.
Here, $\vrho_{\delta,h}\vc{u}_{\delta,h}\otimes
\vc{u}_{\delta,h}\weak \vrho_{\delta}\vc{u}_{\delta}\otimes
\vc{u}_{\delta}$ follows from setting $g_{h} =
\vrho_{\delta,h}\vc{u}_{\delta,h}$ and $f_{h}= \vc{u}_{h, \delta}$
in Lemma \ref{lemma:dLions}, where $g_{h} \weak
g=\vrho_{\delta}\vc{u}_{\delta}$ from the second last result in
\eqref{eq:discconv}.

Next, since
    $\eta_{\delta,h} \in_b L^1(0,T;L^2(\Om))\cap L^1(0,T;W^\frac{3}{2}(\Om))$,
    $\partial_t \eta_{\delta,h} \in_b L^1(0,T;W^{-1,1}(\Om))$ and $W^\frac{3}{2}$ is
    compactly embedded in $L^2$, we can apply Lemma \ref{lemma:simon}
    to conclude that
    \begin{equation}\label{eq:etaconv}
        \eta_{\delta, h} \rightarrow \eta_{\delta} \quad \text{in }L^{2-\epsilon}(0,T;L_{\text{loc}}^2(\Om)),
         \qquad \sqrt{\eta_{\delta,h}} \rightarrow \sqrt{\eta_\delta} \quad \text{in }L^{2}(0,T;L^2_{\text{loc}}(\Om)),
    \end{equation}
Now, by a straight forward application of the H\"older inequality,
we deduce $\Grad \eta_{\delta, h} \in_b L^2(0,T;L^1(\Om))\cap
L^1(0,T;L^\frac{3}{2}(\Om))$. The standard interpolation inequality
provides an estimate of the form $\Grad \eta_{\delta,h}
\in_b L^{\alpha_1}(0,T;L^{\alpha_2}(\Om))$. Thus, in particular
\begin{equation}\label{eq:thetrick}
    \Grad \eta_{\delta,h} \weak \Grad \eta_\delta,
\end{equation}
in the sense of distributions on $(0,T)\times \Om$.

Now, equipped with \eqref{eq:discconv}, \eqref{eq:thetrick}, and
the bounds \eqref{eq:spaces1}, there is no problem with taking the
limit $h \rightarrow 0$ in \eqref{eq:timedisc-vrho} and
\eqref{eq:timedisc-eta} to discover that
\begin{equation}\label{eq:cont-vrho}
-\int_{0}^T\int_{\Om}\vrho_{\delta} (\psi_{t} -
\vc{u}_{\delta}\Grad \psi)\ dxdt = \int_{\Om}\vrho_{\delta,
0}\psi(0,x) \ dx,
\end{equation}
for all $\psi \in C^\infty_{c}([0,T) \times \overline{\Om}))$ and
\begin{equation}\label{eq:cont-eta}
-\int_{0}^T\int_{\Om}\eta_{\delta}\left(\psi_{t} +
(\vc{u}_{\delta} - \Grad\Phi)\Grad\psi\right) - \Grad
\eta_{\delta} \Grad \psi\ dxdt = \int_{\Om}\eta_{\delta,
0}\psi(0,x) \ dx
\end{equation}
for all $\psi \in C^\infty_{c}([0,T) \times \overline{\Om}))$.

Hence, the limiting functions satisfy both the continuity
equation and particle equation in the sense of Definition
\ref{def:level1}.

Similarly, we can go to the limit $h \rightarrow 0$ in
\eqref{eq:timedisc-momentum} to discover in the limit
\begin{equation*}
\begin{split}
&\int_{0}^T\int_{\Om}\vrho_{\delta}\vc{u}_{\delta}\vc{v}_{t} + \vrho_{\delta}\vc{u}_{\delta} \otimes \vc{u}_{\delta}: \Grad \vc{v}_{\delta} + \eta_{\delta}\Div \vc{v} \ dxdt \\
&= \int_{0}^T\int_{\Om}\mu \Grad\vc{u}_{\delta}\Grad \vc{v} +
\lambda \Div \vc{u}_{\delta}\Div \vc{v}
 + (\beta \vrho_{\delta} + \eta_{\delta})\Grad \Phi \vc{v}\ dxdt  \\
&\qquad -  \lim_{h \rightarrow 0} \int_{0}^T\int_{\Om}
p_{\delta}(\vrho_{\delta,h})\Div \vc{v} \ dxdt  -
\int_{\Om}\vc{m}_{0}\vc{v}(0,x) \ dx,
\end{split}
\end{equation*}
for all $\vc{v} \in C^\infty_{c}([0,T) \times \Om)$. Thus, in
order to conclude existence of an artificial pressure solution in
the sense of Definition \ref{def:level1}, we must prove that in
fact
\begin{equation}\label{eq:basicprob1}
\lim_{h \rightarrow 0} \int_{0}^T\int_{\Om}
p_{\delta}(\vrho_{\delta,h})\Div \vc{v} \ dxdt  =
\int_{0}^T\int_{\Om}p_{\delta}(\vrho_{\delta}) \Div \vc{v}\ dxdt.
\end{equation}

\begin{Lemma}
Fix any $\delta > 0$ and let $\{\vrho_{\delta, h}, \vc{u}_{\delta,
h}, \eta_{\delta, h}\}_{h>0}$ be a sequence of time discretized
solutions constructed according to Definition \ref{def:level2} and
\eqref{eq:level2-ext}. Then, there exists a triple
$(\vrho_{\delta}, \vc{u}_{\delta}, \eta_{\delta})$ such that as $h
\rightarrow 0$, $\vrho_{\delta,h} \weak \vrho_{\delta}$ in
$L^\infty(0,T;L^\gamma (\Om)\cap L^6(\Om))$, $\vrho_{\delta,h}
\rightarrow \vrho_{\delta}$ a.e in $\Dom$, $\vc{u}_{\delta,h}
\weak \vc{u}_{\delta}$ in $L^2(0,T;W^{1,2}_{0}(\Om))$,
$\vrho_{\delta, h}\vc{u}_{\delta, h} \weak
\vrho_{\delta}\vc{u}_{\delta}$ in the sense of distributions on
$\Dom$, $\eta_{\delta, h} \weak \eta_{\delta}$ in
$L^1(0,T;W^{1,\frac{3}{2}}(\Om))$, and $\eta_{\delta, h}
\rightarrow \eta_{\delta}$ a.e in $\Dom$, where $(\vrho_{\delta},
\vc{u}_{\delta}, \eta_{\delta})$ is a weak solution to the
artificial pressure approximation scheme in the sense of
Definition \ref{def:level1}.
\end{Lemma}

\bProof In view of the high integrability, and strong convergence
properties, of $\eta_{\delta, h}$, \eqref{eq:basicprob1} can be
proved by the same arguments as those leading to Theorem 7.2 in
\cite{LI4}. Some minor modifications are needed to treat the
unbounded potential $\Phi$. The arguments needed are identical to
those of Section 4.5.2 and will not be given here. From this and
the previous results of this section we can conclude existence of
an artificial pressure solution. It remains to prove the energy
inequality \eqref{i7}.

We start with the following calculation
\begin{equation}\label{eq:multiplying}
    \begin{split}
    &\lim_{h \rightarrow 0}\int_0^t \int_\Om 4|\Grad \sqrt{\eta_{\delta,h}}|^2 +
                4 \Grad \eta_{\delta,h} \Grad \Phi + \eta_{\delta,h}|\Grad \Phi|^2~dxdt \\
    &\qquad \qquad
        = \int_0^t \int_\Om  \eta_{\delta}|\Grad \Phi|^2
                + 4\Grad \eta_{\delta} \Grad \Phi + 4\overline{|\Grad \sqrt{\eta_{\delta}}|^2}~dxdt \\
    &\qquad \qquad \geq \int_0^t \int_\Om \left|2\Grad \sqrt{\eta_{\delta}} +
                     \sqrt{\eta_{\delta}}\Grad \Phi\right|^2~dxdt,
    \end{split}
\end{equation}
where we have used Lemma \ref{lem:prelim} and \eqref{eq:thetrick}. By
taking the limit  $h \rightarrow 0$ in \eqref{eq:disc-energybound}
(using convexity of $z \mapsto z^\gamma$ and $z \mapsto z \log
z$), we obtain

\begin{equation}\label{eq:delta-energy}
    \begin{split}
    &E(\vrho_{\delta}, \vc{u}_{\delta}, \eta_{\delta})(t) \\
    &\qquad
        + \int_{0}^t \int_{\Om}\mu |\Grad \vc{u}_{\delta}|^2 + \lambda |\Div \vc{u}_{\delta}|^2\ dxdt \\
    &\qquad
        + \int_{0}^t\int_{\Om} \left|2\Grad \sqrt{\eta_{\delta}} +
                         \sqrt{\eta_{\delta}}\Grad \Phi\right|^2~dxdt\\
    &\leq E(\vrho_{\delta,0},\vc{u}_{\delta,0}, \eta_{\delta,0}),
    \end{split}
\end{equation}
for any $t \in (0,T)$. \qed

\subsection{Vanishing artificial pressure limit ($\delta \rightarrow 0$)}\label{S3.5}
In the previous subsection we proved that, for each fixed $\delta
> 0$, there exists an artificial pressure solution in the sense of
Definition \ref{def:level1}. Throughout this section, we let
$\Set{\vrho_{\delta}, \vc{u}_{\delta}, \eta_{\delta}}_{\delta >
0}$ be a sequence of such solutions. The aim is now to prove that
this sequence converges as $\delta \rightarrow 0$ to a weak
solution of the fluid--particle interaction model
\eqref{i1}--\eqref{i3} in the sense of Definition \ref{DD2a}. This
will then conclude
 the proof of Lemma \ref{lemma:existence}.

\subsubsection{Energy bounds}
The energy inequality \eqref{i7}, together with Sobolev embedding
and Corollary \ref{corentropybound}, allow us to conclude the
following $\delta$--independent bounds
\begin{equation}\label{eq:spaces2}
\begin{split}
\vrho_{\delta}|\vc{u}_{\delta}|^2 &\in_{b} L^\infty(0,T;L^1(\Om))\cap L^2(0,T;L^{m_{1}}(\Om)), m_{1}>1,\\
\vc{u}_{\delta} &\in_{b} L^2(0,T;W^{1,2}_{0}(\Om)), \\
\vrho_{\delta} &\in_{b} L^\infty(0,T;L^\gamma(\Om)), \\
\vrho_{\delta}\vc{u}_{\delta} &\in_{b} L^\infty(0,T;L^\frac{2\gamma}{\gamma +1}(\Om)), \\
\eta_{\delta}\log \eta_{\delta} & \in_{b} L^\infty(0,T;L^1(\Om)).
\end{split}
\end{equation}
By straight forward applications of the H\"older inequality, using
\eqref{eq:spaces2} and Sobolev embedding, we deduce the bounds
\begin{equation*}
\vr_\delta \vu_\delta \in_b L^2(0,T;L^{m_2}(\Om)), \quad
\vr_\delta \vu_\delta \otimes \vu_\delta \in_b
L^2(0,T;L^{c_2}(\Om)),
\end{equation*}
where
$$
m_2 = \frac{6\gamma}{6 + \gamma},\qquad c_2 = \frac{3\gamma }{3 +
\gamma}.
$$
From \eqref{eq:etaonceandforall}, we also have that
\begin{equation*}
\eta_{\delta} \in_{b} L^2(0,T;L^3(\Om))\cap
L^1(0,T;W^\frac{3}{2}(\Om))
\end{equation*}

Using the above $\delta-$independent estimates, we easily deduce
the weak time control bounds
\begin{equation*}
    \begin{split}
        \partial_{t}\eta_{\delta} &\in_{b} L^1(0,T;W^{-1,\frac{3}{2}}(\Om)), \\
        \partial_{t}\vrho_{\delta} &\in_{b} L^\infty(0,T;W^{-1,\frac{2\gamma}{\gamma+1}}(\Om)), \\
        \partial_{t}(\vrho_{\delta}\vc{u}_{\delta}) &\in_{b} L^1(0,T;W^{-1,1}(\Om)).
    \end{split}
\end{equation*}

\subsubsection{Convergence}
The bounds in the previous subsection asserts the existence of
functions $(\vr, \vu, \eta)$ such that, passing to a subsequence
if necessary,
\begin{align*}
    \vr_{\delta} &\overset{*}{\weak} \vr_\delta \quad \text{in } L^\infty(0,T;L^\gamma(\Om)), \\
    \vu_{\delta} &\weak \vu_\delta \quad \text{in } L^2(0,T;\vc{W}^{1,2}_0(\Om)), \\
    \eta_{\delta} &\weak \eta_\delta \quad \text{in } L^2(0,T;L^3(\Om)).
\end{align*}
As in the previous subsection, Lemma \ref{lemma:dLions} can be
applied to conclude
\begin{equation}\label{eq:contconv}
\vrho_{\delta}\vc{u}_{\delta}\weak \vrho\vc{u}, \qquad
\vrho_{\delta}\vc{u}_{\delta}\otimes \vc{u}_{\delta} \weak
\vrho\vc{u}\otimes \vc{u},
\end{equation}
in the sense of distributions on $\Dom$.

By similar arguments as those leading to \eqref{eq:etaconv}, we
deduce
    \begin{equation*}
        \eta_{\delta} \rightarrow \eta \quad \text{in }L^{2-\epsilon}(0,T;L^2_{\text{loc}}(\Om)),
         \qquad \sqrt{\eta_{\delta}} \rightarrow \sqrt{\eta} \quad \text{in }L^{2}(0,T;L^2_{\text{loc}}(\Om)).
    \end{equation*}
and
\begin{equation}\label{eq:thetrick2}
    \Grad \eta_{\delta} \weak \Grad \eta,
\end{equation}
in the sense of distributions on $(0,T) \times \Om$.

Using \eqref{eq:contconv}, the bounds \eqref{eq:spaces2},
\eqref{eq:thetrick2}, and strong convergence of the initial
conditions, we can take the limit $\delta \rightarrow 0$ in
\eqref{eq:cont-vrho} and \eqref{eq:cont-eta} to discover that
\begin{equation*}
-\int_{0}^T\int_{\Om}\vrho (\psi_{t} - \vc{u}\Grad \psi)\ dxdt =
\int_{\Om}\vrho_{0}\psi(0,x) \ dx,
\end{equation*}
for all $\psi \in C^\infty_{c}([0,T) \times \overline{\Om}))$ and
\begin{equation*}
-\int_{0}^T\int_{\Om}\eta\left(\psi_{t} + (\vc{u} -
\Grad\Phi)\Grad\psi\right) - \Grad \eta\Grad \psi\ dxdt =
\int_{\Om}\eta_{0}\psi(0,x) \ dx
\end{equation*}
for all $\psi \in C^\infty_{c}([0,T) \times \overline{\Om}))$.
Hence, the limiting functions satisfy both the continuity
equation and particle equation in the sense of Definition
\ref{DD2a}.

Similarly, we can go to the limit $\delta \rightarrow 0$ in
\eqref{eq:timedisc-momentum} to discover in the limit
\begin{equation*}
\begin{split}
&\int_{0}^T\int_{\Om}\vrho\vc{u}\vc{v}_{t} + \vrho\vc{u} \otimes \vc{u}: \Grad \vc{v}_{\delta} + \eta\Div \vc{v} \ dxdt \\
&= \int_{0}^T\int_{\Om}\mu \Grad\vc{u}\Grad \vc{v} + \lambda \Div \vc{u}\Div \vc{v} + (\beta \vrho + \eta)\Grad \Phi \vc{v}\ dxdt  \\
&\qquad -  \lim_{\delta \rightarrow 0} \int_{0}^T\int_{\Om}
p(\vrho_{\delta})\Div \vc{v} \ dxdt  -
\int_{\Om}\vc{m}_{0}\vc{v}(0,x) \ dx,
\end{split}
\end{equation*}
for all $\vc{v} \in C^\infty_{c}([0,T) \times \Om)$. Thus, in
order to conclude existence of a weak solution in the sense of
Definition \ref{DD2a}, it remains to prove that
\begin{equation*}
\lim_{\delta \rightarrow 0} \int_{0}^T\int_{\Om}
p(\vrho_{\delta})\Div \vc{v} \ dxdt  =
\int_{0}^T\int_{\Om}p(\vrho) \Div \vc{v}\ dxdt,
\end{equation*}
together with the energy inequality \eqref{i7}. Consequently, we
are faced with a similar situation as in the previous subsection.
The main difference is that  we now only have  $\gamma >
\frac{N}{2}$.
Since $\eta_{\delta}$ enjoys both high integrability and
compactness, the proof follows by a small extension of Feireisl's
arguments
 in \cite{FNP}.
First, we establish higher integrability of the density on the
entire domain $\Om$.

\begin{Lemma}
Let  $\Set{\vrho_{\delta},\vc{u}_{\delta}, \eta_{\delta}}_{\delta
> 0}$ be a sequence of artificial pressure solutions in the sense
of Definition \ref{def:level1}. Then, there exists a constant
$c(T)$, independent of $\delta$, such that
\begin{equation*}
\int_0^T\int_\Om \vrho_{\delta}^{\gamma + \theta}~dxdt \leq c(T),
\end{equation*}
where $\Theta = \min\{\frac{2}{3}\gamma - 1, \frac{1}{4}\}$.
\end{Lemma}

\bProof Since $\eta_{\delta} \in_{b} L^2(0,T;L^3(\Om))\cap
L^1(0,T;W^{1,\frac{3}{2}}(\Om))$, the addition of $\eta_{\delta}$
in the equations does not impose any potential problems.
Consequently, if $\Om$ is bounded, the proof follows by the
classical arguments (\cite{LI4,FNP}).

If $\Om$ is unbounded,  then $\Grad \Phi$ is no longer integrable
and we cannot simply apply existing results. To prove the bound in
this case, let $\Delta^{-1}$ be the inverse Laplacian realized
using Fourier multipliers (see \cite{LI4, EF70} for details). For
each fixed $\delta >0 $, let the test-function $\vc{v}_\delta$ be
given as
\begin{equation*}
    \vc{v}_\delta = \Grad \Delta^{-1}\vr^\theta_\delta.
\end{equation*}
By the requirements on $\theta$, we have in particular
$\vr_{\delta}^\theta \in_b L^\infty(0,T;L^s(\Om))$,
$$
    s = \max\left\{\frac{3\gamma }{2\gamma - 3}, 4\gamma\right\}.
$$
Thus, $ \vc{v}_\delta \in_b L^\infty(0,T;W^{1,s}(\Om)) \cap
L^\infty(0,T;L^\infty(\Om)). $

Next, since $(\vr_\delta, \vu_\delta)$ is a renormalized solution
to the continuity equations, \eqref{ii1} with $B(\vr_\delta) =
\vr_\delta^\theta$ states
\begin{equation*}
    \partial_t \vr_\delta^\theta = -\Div (\vr_\delta^\theta \vu) - (\theta-1)\vr_\delta^\theta \Div \vu,
\end{equation*}
in the sense of distributions on $(0,T)\times \Om$.  For
notational convenience, we observe that
\begin{equation*}
    \begin{split}
        \|\partial_t \vc{v}_\delta\|_{L^p(0,T;L^q(\Om))}
        &= \|\Grad \Delta^{-1} \partial_t \vr^\theta_\delta\|_{L^p(0,T;L^q(\Om))} \\
        &\leq \|\vr_{\delta}^\theta \vu_\delta\|_{L^p(0,T;L^q(\Om))} +  \|\vr^\theta \Div \vu\|_{L^p(0,T;L^{r}(\Om))},
    \end{split}
\end{equation*}
for appropriate $1\leq p,q \leq \infty$ and $r^* = q$.

Next, we apply $\vc{v}_\delta$ as test function for the momentum
equation to obtain
\begin{equation*}
    \begin{split}
    \int_0^T\int_\Om a\vr_\delta^{\gamma + \theta}~dxdt
    &=  -\int_0^T\int_\Om (\vr_\delta \vu_\delta) \partial_t \vc{v}_\delta + \vr_\delta \vu_\delta \otimes \vu_\delta :\Grad \vc{v}_\delta~dxdt \\
    &\quad + \int_0^T\int_\Om \mu \Grad \vu_\delta \Grad \vc{v}_\delta + \lambda \Div \vu_\delta \Div \vc{v}_\delta~dxdt \\
    &\quad- \int_0^T\int_\Om \eta_\delta \vr_\delta^\theta  - (\vr_\delta\beta + \eta_\delta)\Grad \Phi_\delta \vc{v}_\delta~dxdt
     - \int_\Om \vc{m}_0\vc{v}_\delta(0,\cdot)~dx \\
    &:= I_1 + I_2 + I_3.
    \end{split}
\end{equation*}
To conclude it remains to bound $I_1$, $I_2$, and $I_3$,
independently of $\delta$. We start with the $I_1$ term:
\begin{equation*}
    \begin{split}
        |I_1| &:= \left|\int_0^T\int_\Om \vr_\delta \partial_t \vc{v} + \vr_\delta \vu_\delta \otimes \vu_\delta :\Grad \vc{v}_\delta~dxdt\right|\\
            &\leq
            \|\vr_\delta \vu_\delta \|_{L^\infty(0,T:L^\frac{2\gamma}{\gamma+1}(\Om))}
                \|\vr_{\delta}^\theta \vu_\delta\|_{L^1(0,T;L^\frac{2\gamma}{\gamma-1}(\Om))} \\
            & + \|\vr_\delta \vu_\delta \|_{L^2(0,T:L^{m_2}(\Om))}
                C(T)\|\vr^\theta \Div \vu\|_{L^2(0,T;L^{r}(\Om))} \\
            &\qquad + \|\vr_\delta \vu_\delta \otimes \vu_\delta\|_{L^2(0,T;L^{c_2}(\Om))}\|\Grad \vc{v}_\delta\|_{L^2(0,T;L^{c_2'}(\Om))},
    \end{split}
\end{equation*}
where
$$
r= \frac{6\gamma }{7\gamma -6 },\qquad r^* = (m_2)', \qquad  c_2'
= \frac{3\gamma }{2\gamma - 3} \leq s.
$$
Now, we estimate
\begin{align*}
    \|\vr^\theta_\delta \vu_\delta\|_{L^1(0,T;L^\frac{2\gamma}{\gamma-1}(\Om))}
        &\leq C(T)\|\vr_\delta^\theta\|_{L^\infty(0,T;L^s(\Om))}\|\vu_\delta\|_{L^2(0,T;L^{6}(\Om)}, \\
    \|\vr_\delta^\theta \Div \vu_\delta\|_{L^2(0,T;L^{r}(\Om))}
        &\leq  C\|\Div \vu_\delta\|_{L^2(0,T;L^2(\Om))}\|\vr^\theta_\delta\|_{L^\infty(0,T;L^s(\Om))},
\end{align*}
and hence conclude that
$
|I_1| \leq C(T).
$
Next, we easily deduce the bound
$
|I_2| \leq C(T),
$
and it only remains to bound $I_3$.
\begin{equation*}
    \begin{split}
        |I_3| &\leq \|\eta_\delta\|_{L^2(0,T;L^3(\Om))}\|\vr^\theta\|_{L^2(0,T;L^\frac{3}{2}(\Om))} \\
                & \qquad + \left\|(\beta \vr_\delta + \eta_\delta) \Grad \Phi \right\|_{L^1(0,T;L^1(\Om))}\|\vc{v}_\delta\|_{L^\infty(0,T;L^\infty(\Om))}
                    + \|\vc{m}_0\|_{L^1(\Om)}\|\vc{v}_\delta\|_{L^\infty(\Om)}\\
              &\leq C(T)\left( 1 + \left\|(\beta \vr_\delta + \eta_\delta) \Grad \Phi \right\|_{L^\infty(0,T;L^1(\Om))}\right).
    \end{split}
\end{equation*}
Using the energy estimate \eqref{eq:delta-energy}, Corollary
\ref{corentropybound}, and the requirements \eqref{potentialreq}
on the potential, we readily deduce
\begin{equation*}
    \begin{split}
    &\sup_{t \in (0,T)}\int_\Om |\beta \vr_\delta + \eta_\delta||\Grad \Phi|~dx \\
    & \qquad   \leq \|\Grad \Phi\|_{L^\infty(B(0,R))}\sup_{t \in (0,T)}\int_{B(0,R)}\beta \vr_\delta + \eta_\delta~dx \\
    &\qquad \qquad + C\sup_{t \in (0,T)}\int_{\Om \backslash B(0,R)}(\beta \vr_\delta + \eta_\delta)\Phi~dxdt
    \leq C(T),
    \end{split}
\end{equation*}
which brings the proof to an end.

\qed

The following lemma concludes the proof of Lemma
\ref{lemma:existence}.

\begin{Lemma}
Let  $\Set{\vrho_{\delta},\vc{u}_{\delta}, \eta_{\delta}}_{\delta
> 0}$ be a sequence of artificial pressure solutions in the sense
of Definition \ref{def:level1}. Then, there exists a triple
$(\vrho, \vc{u}, \eta)$ such that as $\delta \rightarrow 0$,
$\vrho_{\delta} \weak \vrho$ in $L^\infty(0,T;L^\gamma (\Om))$,
$\vrho_{\delta} \rightarrow \vrho$ a.e in $\Dom$, $\vc{u}_{\delta}
\weak \vc{u}$ in $L^2(0,T;W^{1,2}_{0}(\Om))$,
$\vrho_{\delta}\vc{u}_{\delta} \weak \vrho\vc{u}$ in the sense of
distributions on $\Dom$, $\eta_{\delta} \weak \eta$ in
$L^1(0,T;W^{1,\frac{3}{2}}(\Om))$, and $\eta_{\delta} \rightarrow
\eta$ a.e in $\Dom$, where $(\vrho, \vc{u}, \eta)$ is a weak
solution to the particle interaction model in the sense of
Definition \ref{DD2a}.

Moreover, the total fluid mass and particle mass given by
$$
M_\vr(t) = \int_\Om \vr(t, \cdot)~dx \qquad \mbox{and} \qquad
M_\eta (t) = \int_\Om \eta(t, \cdot)~dx,
$$
respectively, are constants of motion.
\end{Lemma}
\bProof
    1. Since $\eta_{\delta} \rightarrow \eta$ a.e in $\Dom$, the presence of the $\eta_{\delta}$
    variable does not impose any extra difficulties as compared with the corresponding situation for the
    barotropic compressible
    Navier--Stokes equations.
    Consequently, strong convergence of the density and  existence of a solution follows as
    in \cite[Section 4.2]{FNP}. In the case $\Om$ bounded, this concludes the proof.
    For an unbounded domain, the potential is non-integrable and
    it remains to prove the energy estimate \eqref{i7}.

    Since the mappings $z \mapsto z^\gamma$ and $z \mapsto z \log z$ are
    convex, we can apply Lemma \ref{lem:prelim} to conclude
    that, for a.e $t \in (0,T)$,
    \begin{equation}\label{eq:lp1}
        E(\vr, \vu, \eta) \leq \lim_{\delta \rightarrow 0}E(\vr_\delta, \vu_\delta, \eta_\delta).
    \end{equation}
    Next, similar to \eqref{eq:multiplying}, we deduce,
    \begin{equation}\label{eq:lp2}
        \begin{split}
            & \int_0^T \int_K \left|2\Grad \sqrt{\eta} + \sqrt{\eta}\Grad \Phi  \right|^2\ dxdt  \\
            &\qquad \leq \lim_{\delta \rightarrow  0}\left[
                \int_0^T \int_K 4|\Grad \sqrt{\eta_\delta}|^2 + 4\Grad \eta_\delta \Grad \Phi + \eta_\delta |\Grad \Phi|^2~dxdt\right] \\
            &\qquad = \lim_{\delta \rightarrow  0}\left[\int_0^T\int_\Om\left|2\Grad \sqrt{\eta_\delta} + \sqrt{\eta_\delta}\Grad \Phi  \right|^2 \ dxdt\right] \\
            &\qquad \leq E(\vr_0, \vu_0, \eta_0).
        \end{split}
    \end{equation}
    for any compact $K$.

    Finally, using \eqref{eq:lp1} and \eqref{eq:lp2}, we can
    send $\delta \rightarrow 0$ in \eqref{eq:delta-energy} to obtain
    \begin{equation*}
        \begin{split}
            &E(\vrho, \vc{u}, \eta)(t) \\
            &\qquad
                + \int_{0}^t \int_{\Om}\mu |\Grad \vc{u}|^2 + \lambda |\Div \vc{u}|^2\ dxdt \\
            &\qquad
                + \int_{0}^t\int_{\Om}\eta\left|\Grad \Phi + \Grad \log \eta \right|^2\ dxdt\\
            &\leq E(\vrho_{0},\vc{u}_{0}, \eta_{0}).
        \end{split}
    \end{equation*}

 2. Let us now prove that $M_\vr$ and $M_\eta$ are constants of motion. To prove this, we make use of Lemma \ref{lemcontmass}
 and conclude that for any $\epsilon > 0$ there exists $R > 0$ such that
\begin{equation*}
    \int_\Om \eta_{\delta, h}~dx = \int_{\Om \cap B(0,R)} \eta_{\delta, h}~dx + \mathcal{O}(\epsilon).
\end{equation*}
By sending $\delta, h \rightarrow 0$, we gather
\begin{equation*}
    M_\vr(0) = \int_\Om \eta_0~dx = \lim_{\delta, h \rightarrow 0}\int_\Om \eta_{\delta, h}~dx
    = \int_{\Om \cap B(0,R)} \eta~dx + \mathcal{O}(\epsilon),
\end{equation*}
which implies that $M_\eta$ is a constant of motion.

From the energy estimate \eqref{eq:disc-energybound} and Corollary
\ref{corentropybound} it follows that (for a.e $t \in (0,T)$)
\begin{equation*}
    \int_\Om \vr_{h,\delta}\Phi~dx \leq C.
\end{equation*}
By the requirements of the potential \eqref{potentialreq}, we can
for any $M > 0$ determine a radius $R$ such that $\Phi(x) \geq M$,
for all $|x|> R$. Hence,
\begin{equation*}
    \int_\Om \vr_{\delta, h}~dx = \int_{ B(0,R)} \vr_{\delta, h}~dx + \int_{\Om \setminus B(0,R)} \vr_{\delta, h}~dx,
\end{equation*}
where
\begin{equation*}
     \int_{\Om \setminus B(0,R)} \vr_{\delta, h}~dx \leq M^{-1}\int_{\Om \setminus B(0,R)} \Phi \vr_{\delta, h}~dx
    \leq M^{-1}C.
\end{equation*}
Since this holds for any $M>0$, we  conclude that $M_\vr(t)$ is a
constant of motion.

\qed


%

\section{Large-time Asymptotics}\label{S4}
The analysis in the previous section yields
\begin{Lemma}\label{lemma:zero}
Under the hypotheses of Theorem \ref{T2}, we have
$$\!\!\!\!\! \lim_{\! \tau \to \infty\!} \int_{\tau-1}^{\tau+2} \!\!\!\!\!\|\nabla \vu\|^2_{\!L^{p_1}\!(\Omega)} \!+\! \|\vr |\vu|^2\|_{\!L^{p_2}\!(\Omega)} \!+\!
\|\vr |\vu|\|^2_{\!L^{p_3}\!(\Omega)}\! 
+\! \| |\eta |\vu|\|_{\!L^{p_4}\!(\Omega)}dt\! =0,$$
\begin{equation}\label{eq:entropyprod}
    \lim_{\tau \to \infty}\int_{\tau-1}^{\tau+2}\int_\Om |2\Grad \sqrt{\eta} + \sqrt{\eta}\Grad \Phi|^2~dxdt = 0,
\end{equation}
and
\begin{equation*}
    \int_{\tau-1}^{\tau+2}\int_\Om \vr^{\gamma + \theta}~dxdt \leq C,\quad \text{for all }\tau > 1.
\end{equation*}
Here,
$$p_1 =2, \quad p_2 = \frac{3\gamma}{\gamma +3},\quad  p_3 = \frac{6 \gamma}{\gamma + 6},
\quad  p_4 = 4,$$
$$
\theta = \min\Big\{\frac{2}{3}\gamma-1,\frac{1}{4}\Big\}.
$$
\end{Lemma}

Let the sequence $\{t_n\}_{n=1}^\infty$ be such that
$t_n\rightarrow\infty$ and define the sequences:
$$
\vu_n(t,x)=\vu(t+t_n,x),~~\vr_n(t,x)=\vr(t+t_n,x),~~\eta_n(t,x)=\eta(t+t_n,x),
\,\,\, t\in(-1,2), \, x \in \Omega.
$$
Observe that, for each fixed $n$, the triple $(\vr_n, \vu_n,
\eta_n)$ is a weak solution to particle interaction model in the
sense of Definition \ref{DD2a}.

By virtue of the previous lemma, and bounded energy, we can
conclude the existence of functions $\vr_s$, $\overline{p}$, and
$\eta_s$, such that
\begin{eqnarray}
& & \vr_n = \vr(t + \tau_n) \to \vr_s \tab  weakly\,\,  in \,\, L^{\gamma}((-1,2) \times \Omega)  \nonumber\\
& & p_n = p(t + \tau_n) \to \overline{p} \tab  weakly\,\,  in \,\, L^{1}((-1,2) \times \Omega)  \label{eq:conv1}\\
& & \eta_n = \eta(t + \tau_n) \to \eta_s \tab weakly\,\,  in \,\,
L^2((-1,2), L^3(\Omega))\nonumber.
\end{eqnarray}


\begin{Lemma}
Given \eqref{eq:conv1} and the hypothesis of Theorem \ref{T2},
$\eta_n \rightarrow \eta_s$ a.e in $(-1,2)\times \Om$, where
$\eta_s$ is given by
\begin{equation*}
    \eta_s = C_\eta\exp(-\Phi), \quad C_\eta = \frac{\int_\Om \exp(-\Phi)~dx}{\int_\Om \eta_0~dx},
\end{equation*}
which is the unique solution to the problem
\begin{equation}\label{eq:linear}
    \Grad \eta_s = - \eta_s\Grad \Phi, \quad \int_\Om \eta_s~dx = \int_\Om \eta_0~dx.
\end{equation}

\end{Lemma}
\bProof From the particle equation \eqref{i3}, we observe that
\begin{equation*}
    \partial_t\eta_n = -\Div(\eta_n \vu_n) + \Div\Big[\sqrt{\eta_n}\left(\sqrt{\eta_n}\Grad \Phi + 2\Grad \sqrt{\eta}\right)\Big]
\end{equation*}
in the sense of distributions on $(-1,2)\times \Om$. Hence,
\begin{equation*}
    \begin{split}
        \int_{-1}^2\|\partial_t \eta_n\|_{W^{-1,1}(\Om)}~dt
         \leq& \|\eta_n\vu_n\|_{L^1(-1,2;L^4(\Om))} \\
        & + \|\eta_n\|_{L^2(-1,2;L^2(\Om))}\|\sqrt{\eta_n}\Grad \Phi + 2\Grad \sqrt{\eta}\|_{L^2(-1,2;L^2(\Om))}.
    \end{split}
\end{equation*}
Lemma \ref{lemma:zero} tells us that the right-hand side converges
to zero as $n \rightarrow \infty$. Since in addition $\eta_n \in_b
L^1(-1,2;W^\frac{3}{2}(\Om))$, Lemma \ref{lemma:simon} can be
applied to conclude that $\eta_{n} \rightarrow \eta_s$ a.e in
$(-1,2) \times \Om$. In addition, as in \eqref{eq:thetrick}, we
find that $\Grad \eta_n \weak \Grad \eta$ in the sense of
distributions on $(-1,2)\times \Om$.

Then, as in \eqref{eq:lp2}, we calculate
\begin{equation*}
    \begin{split}
        0 &\leq \int_{-1}^2 \int_K \left| \Grad \sqrt{\eta_s} + \sqrt{\eta_s}\Grad \Phi\right|^2~dxdt \\
        & \leq \lim_{n \rightarrow \infty}\left[\int_{-1}^{2}\int_K |2\Grad \sqrt{\eta_n} + \sqrt{\eta_n}\Grad \Phi|^2~dxdt\right] \leq 0,
    \end{split}
\end{equation*}
for any compact $K$. The last equality follows directly from
\eqref{eq:entropyprod}. Hence, we can conclude
\begin{equation*}
    \Grad \sqrt{\eta_s} = -\sqrt{\eta_s}\Grad \Phi \quad \text{a.e in }(0,T)\times \Om.
\end{equation*}
Clearly, this means that also
\begin{equation}\label{eq:ae}
    \Grad \eta_s = -\eta_s \Grad \Phi\quad \text{a.e in }(0,T)\times \Om.
\end{equation}
Next, we let $n \rightarrow \infty$ in \eqref{i3} and apply Lemma
\ref{lemma:zero} and \eqref{eq:ae} to obtain
\begin{equation*}
    \int_{-1}^2 \int_\Om \eta_s \phi_t~dxdt = 0, \quad \forall \phi \in C_c^\infty((-1,2)\times \Om).
\end{equation*}
Hence, $\eta_s$ is \emph{independent of time}. Using this, we
deduce
\begin{equation*}
    \begin{split}
        \int_\Om \eta_0~dx
            &= \lim_{n \rightarrow \infty}\frac{1}{3}\int_{-1}^2 \int_\Om \eta_n ~dxdt\\
            &= \lim_{n \rightarrow \infty}\left[\frac{1}{3}\int_{-1}^2 \int_{B(0,R)} \eta_n ~dxdt
                + \frac{1}{3}\int_{-1}^2 \int_{\Om \setminus B(0,R)}\eta_n~dxdt\right] \\
            &= \int_{B(0,R)} \eta_s~dx + \lim_{n \rightarrow \infty}\left[\frac{1}{3}\int_{-1}^2 \int_{\Om \setminus B(0,R)}\eta_n~dxdt\right]
    \end{split}
\end{equation*}
for any ball $B(0,R)$ of finite radius $R$. Now, in view of
\eqref{potentialreq}, for any $M>0$ there exists a radius $R$ such
that $\Phi \geq M$. Consequently,
\begin{equation*}
    \frac{1}{3}\int_{-1}^2 \int_{\Om \setminus B(0,R)}\eta_n~dxdt
    \leq \frac{1}{M}\int_{\Om}\eta_n \Phi~~dxdt \leq M^{-1}D,
\end{equation*}
where $D$ is the constant appearing in Corollary
\ref{corentropybound}. Letting $M \rightarrow \infty$, we conclude
\begin{equation}\label{eq:etacons}
 \int_\Om \eta_s~dx = \int_\Om \eta_0~dx.
\end{equation}
Since \eqref{eq:linear} is linear in $\eta_s$ and $\Grad \eta_s$,
this concludes the proof. \qed

By taking the limit $n \rightarrow \infty$ in the continuity
equation \eqref{i1}, keeping in mind Lemma \ref{lemma:zero}, we
obtain
\begin{equation*}
    \int_{-1}^2\int_\Om \vr_s \phi_t~dxdt = 0, \quad \forall \phi \in C_c^\infty((-1,2)\times \Om).
\end{equation*}
Thus, $\vr_s$ is \emph{independent of time}. By virtue of
Corollary \ref{corentropybound} and the energy estimate
\eqref{ii4}, we have that $\beta \vr_s\Phi \in
L^\infty(0,T;L^1(\Om))$. Then, we can redo the arguments leading
to \eqref{eq:etacons}, with $\beta \vr$ replacing $\eta$, to
obtain
\begin{equation}\label{eq:masscons}
    \int_\Om \vr_s~dx = \int_\Om \vr_0~dx.
\end{equation}
Note that we do not lose mass in the large-time limit
(cf~\cite{FP9}).
\begin{Lemma}
Given the convergences \eqref{eq:conv1}, $\vr_n \rightarrow \vr_s$
in $L^\gamma((-1,2)\times \Om)$, where
\begin{equation}\label{eq:theform}
 \vr_s = \left(\frac{\gamma-1}{a\gamma}[-\beta \Phi + C_\vr]^+ \right)^\frac{1}{\gamma-1},
\end{equation}
$C_\vr$ is uniquely determined by $\int_\Om \vr_0~dx$, and $\vr_s$
is the unique solution of
$$
\Grad p(\vr_s) = -\beta\vr_s \Grad \Phi \quad \text{in }\Om, \quad
\int_\Om \vr_s~dx = \int_\Om \vr_0~dx
$$
\end{Lemma}

\bProof With the already obtained bounds, there is no problem with
passing to the limit in the momentum equation \eqref{i2} to obtain
\begin{equation*}
    \Grad \overline{p} + \Grad \eta = -\beta \vr_s \Grad \Phi - \eta_s \Grad \Phi,
\end{equation*}
in the sense of distributions on $\Om$. By virtue of the previous
lemma,
\begin{equation*}
    \Grad \overline{p} = - \beta \vr_s \Grad \Phi.
\end{equation*}
However, then we are in the same situation as in \cite{FP9} and
the arguments of \cite[Proposition 4.1]{FP9} can be applied to
obtain strong convergence $\vr_n \rightarrow \vr_s$ in
$L^\gamma((-1,2)\times \Om)$. Consequently, $\overline{p} =
p(\vr_s)$, and, in view of \eqref{eq:masscons}, $\vr_s$ solves
\begin{equation*}
    \Grad p(\vr_s) = - \beta \vr_s \Grad \Phi,\qquad \int_\Om \vr_s~dx = \int_\Om \vr_0~dx
\end{equation*}
in the sense of distributions. According to \cite[Theorem
2.1]{FPJAVELL}, $\vr_s$ is then of the form \eqref{eq:theform},
where $C_\vr$ is uniquely determined by $\int_\Om \vr_s~dx$. \qed

Lemmas 5.2 and 5.3 concludes the proof of Theorem \ref{T2}.

\section*{Acknowledgments}
This paper was written as part of  the international research
program on Nonlinear Partial Differential Equations at the Centre
for Advanced Study at the Norwegian Academy of Science and Letters
in Oslo during the academic year 2008Ð-09. The authors gratefully
acknowledge the hospitality and support of the Centre for Advanced
Study where this research was performed. Jos\'e A. Carrillo
acknowledges support from the project MTM2008-06349-C03-03 DGI-MCI
(Spain) and 2009-SGR-345 from AGAUR-Generalitat de Catalunya. KT
acknowledges also the support in part by the National Science
Foundation under the awards PECASE DMS 0239063,  DMS 0807815 and
DMS 0405853.

\end{document}